\UseRawInputEncoding
\documentclass[12pt]{article}
\usepackage{bm}
\usepackage{mathptmx}
\usepackage{amsfonts}
\usepackage{latexsym}
\usepackage{graphicx}
\usepackage{color}
\usepackage{ifpdf}
\usepackage[titletoc]{appendix}
\newtheorem{theorem}{Theorem}[section]

\newtheorem{lemma}[theorem]{Lemma}

\newtheorem{definition}{Definition}[section]
\newtheorem{remark}{Remark}[section]
\usepackage{amsfonts,amssymb,eucal,amsmath}
\title{\bf On the Stability of Lagrange Relative Equilibrium  in the Planar Three-body Problem   \thanks{Supported by NSFC(No.11701464)}}

\author{{ Xiang Yu\footnote{Email:xiang.zhiy@foxmail.com, xiang.zhiy@gmail.com}} \\
\small \it Department of
Mathematics, Sichuan University, \\
\small \it Chengdu 610064,  China}

\date{}

\begin{document}
\maketitle

\begin{abstract}
Since  the strong degeneracies present in the N-body problem, even in the basic
case of the planar three-body problem,  nobody inspects the problem of nonlinear stability of Lagrange relative equilibrium. We introduce a new coordinate system to reduce degeneracies according to intrinsic symmetrical characteristic of the $N$-body problem, then we prove  that Lagrange relative equilibrium is stable in the sense of measure, provided it is spectrally stable and except six special resonant cases.  Indeed, under this condition, there are abundant KAM invariant tori or quasi-periodic solutions near Lagrange relative equilibrium. Furthermore, these tori or quasi-periodic solutions  form a set whose relative measure rapidly tends to 1.  We also  prove that Lagrange relative equilibrium  is exponentially stable for almost every choice of  masses in the sense of measure, provided it is spectrally stable; and topologically, this is also right for  a large open subset  of spectrally stable space of masses. \\


{\bf 2010AMS Subject Classification} {70F10 \and 34D20 \and 70H14 \and 70K20 \and 70F15  \and 70K45}.

\end{abstract}

\section{Introduction}
\ \ \ \
We consider $N $ particles with positive mass moving in an Euclidean space ${\mathbb{R}}^2$ interacting under
the laws of universal gravitation. Let the $k$-th particle have mass $m_k$ and position
$\mathbf{r}_k \in {\mathbb{R}}^2$ ($k=1, 2, \cdots, N$), then the equation of motion of the $N$-body problem is written
\begin{equation}\label{eq:Newton's equation1}
m_k \ddot{\mathbf{r}}_k =\sum_{1 \leq j \leq N, j \neq k} \frac{m_km_j(\mathbf{r}_j-\mathbf{r}_k)}{|\mathbf{r}_j-\mathbf{r}_k|^3}, ~~~~~~~~~~~~~~~k=1, 2, \cdots, N.
\end{equation}
where $|\cdot|$  denotes the Euclidean norm in ${\mathbb{R}}^2$. Since these equations
are invariant by translation, we can assume that the center of mass stays at the origin.

The problem of stability of motion of the $N$-body problem is one of  the major problems in mathematics or even in science. However, there are various stability concepts  in modern literature, e.g., Lyapunov stability, Lagrange stability, Poisson stability,  linear stability,  spectral stability, orbital stability, effective stability and stability in the sense of measure (i.e., KAM stability) etc. Beyond doubt, Lyapunov stability is what we all want most and the most difficult to obtain. Note that orbital stability is essentially a kind of stability in the sense of Lyapunov.

For the $N$-body problem, the first result on stability was due to Laplace and Lagrange: 
 the Laplace
theorem on the absence of secular perturbations of the semiaxes. The result was on the stability of the solar system in the sense of Lagrange stability. The method of Laplace is a perturbation analysis neglecting terms in the mass of second-order and higher. Then Poisson and Jacobi 
extended the perturbation analysis  to third-order
terms in the mass, and concluded that Lagrange stability of the solar system is not guaranteed by the method the truncation of the order in the perturbation
analysis.

Next major breakthrough is the well known Arnold's theorem \cite{arnol1963small,fejoz2004demonstration,chierchia2011planetary}, a success of modern celebrated KAM theory, which claims that: for sufficiently small masses
of the planets, Lagrange stability of the solar system is  guaranteed for a set of positive Lebesgue measure of initial conditions. In particular, in the planar restricted circular three-body problem, by the well known Kolmogorov-Arnold's theorem \cite{arnold1963proof}, if the
mass of Jupiter is sufficiently small, then the motion of the the asteroid is stable in the sense of Lagrange stability for most of the initial conditions.

Since the independent work  of Gascheau in 1843 \cite{gascheau1843examen} and Routh in 1875 \cite{routh1874laplace} on linear stability of Lagrange relative equilibrium
of the three-body problem, there are a good deal of work on linear stability of the $N$-body problem,  please see \cite[etc]{moeckel1995linear,roberts1999spectral,roberts2002linear,meyer2005elliptic,hu2010morse,hu2014linear}
and the references therein.

The problem of Lyapunov stability has been solved only for Lagrange relative equilibrium  in the planar restricted circular three-body problem, by a great amount of work based upon KAM theory in the 1970s. Please see \cite[etc]{deprit1967stability,leontovich1962stability,markeev1978libration,sokolskii1978proof,meyer1986stability}
and the references therein. However, due to the possibility of the well known Arnold
diffusion, the method has the limitation that, the number of degrees
of freedom of the problem is not more than 2.

For
physical application, it is natural to consider a sort of  ``effective stability", i.e.,  stability up to finite but long times.  More precisely, for any solution $q(t)$ of a system, with initial condition in a small $\epsilon$-neighbourhood of the equilibrium point $q_0$ of the system,  one could guaranty the  estimate $|q(t)-q_0| = O(\epsilon^a)$ for all times $|t|\leq T(\epsilon)$, where $a$ is some positive number in the
interval $(0,1)$, and $T(\epsilon)$ is a ``large" time such that $T(\epsilon)=O(\frac{1}{\epsilon^b})$ or even more stronger $T(\epsilon)\sim exp(\frac{1}{\epsilon^b})$ for some positive number $b$. The latter stronger form of stability is well known as exponential stability. The exponential stability was first stated by Moser \cite{Moser1955} and Littlewood \cite{littlewood1959equilateral,littlewood1959lagrange} in particular
cases.  A general framework in this direction was developed by Nekhoroshev \cite{nekhorshev1977exponential}. Then a great amount of work focus on the effective stability of Lagrange relative equilibrium  in the planar or spatial restricted circular three-body problem, please see \cite[etc]{Giorgilli1989Effective,celletti1990stability,giorgilli1997stability,benettin1998nekhoroshev,giorgilli1998problem}
and the references therein.
As a matter of fact, the stability over long times has been investigated by Birkhoff \cite{birkhoff1927dynamical} using the
method of normal form going back to Poincar\'{e}.

Although a great amount of work on  stability has been done, however, even in the basic
case of the planar three-body problem, it seems that nobody  consider  the problem of nonlinear stability of Lagrange relative equilibrium. One of the reason may be that one could not reduce the degeneracy of the equations of motion caused by integral of the $N$-body problem well.
Inspired by the work on  the problem of the infinite spin \cite{yu2019problem}, we find that the so-called moving   coordinates introduced in \cite{yu2019problem} are suitable for describing    orbits near central configurations in the planar $N$-body problem, so it is natural to utilize the moving coordinates  to investigate  stability of    relative equilibria. In the moving   coordinates, the degeneracy of the equations of motion according to intrinsic symmetrical characteristic of the $N$-body problem can easily be reduced. In fact,  one can simply reduce the degeneracy corresponding to rotation symmetry of the $N$-body problem, and obtain  practical equations of motion in a small neighbourhood of a relative equilibrium. This is one important reason we can investigate  nonlinear stability of relative equilibria.

The paper is structured as follows. In \textbf{Section 2}, we recall some notations,  and some preliminary results including the moving   coordinates in \cite{yu2019problem}, and give equations of motion by the moving   coordinates.  In \textbf{Section 3}, we simply  discuss  orbital and linear stability of relative equilibria of the planar $N$-body problem, to prepare for investigating KAM stability and effective stability  of Lagrange relative equilibrium  in the planar three-body problem. In \textbf{Section 4}, we recall some necessary classical aspects of Hamiltonian system. In \textbf{Section 5}, we give the  Birkhoff normal
form of the Hamiltonian near Lagrange triangular
point. Although the construction of the normal form is  simple in concept but technically complicated in operations, which require some computer assistance; the analysis was performed with the aid of Mathematica. In \textbf{Section 6}, we investigate KAM stability  of Lagrange relative equilibrium, in particular, it is shown that there are a great quantity of quasi-periodic solutions in a small neighbourhood of Lagrange relative equilibrium. Finally, in \textbf{Section 7}, we investigate  effective stability of Lagrange relative equilibrium.

\section{Preliminaries}
\label{Preliminaries}
\indent\par
In this section we recall some notations and definitions given in \cite{yu2019problem}. In particular, we recall the so-called moving   coordinates introduced to study collision orbits in \cite{yu2019problem}, which would be quite useful for the investigation of the stability of  relative equilibrium solutions. 

Let $(\mathbb{R}^2)^N $ denote the space of configurations for $N$ point particles in Euclidean space $\mathbb{R}^2$: $ (\mathbb{R}^2)^N = \{ \mathbf{r} = (\mathbf{r}_1,  \cdots, \mathbf{r}_N):\mathbf{r}_j \in \mathbb{R}^2, j = 1,  \cdots, N \}$. We remark that, when necessary, one may identify $\mathbb{R}^2$ with $\mathbb{C}$ and $(\mathbb{R}^2)^N$ with $\mathbb{C}^N$ and so on; in particular, $\mathbf{S}^1=\{e^{\mathbf{i}\theta}|\theta\in \mathbb{R}\}$, the unit circle in $\mathbb{C}$, is identified with the special orthogonal group $S\mathbb{O}(2)$ of the plane, where $\mathbf{i}$ is the imaginary unit. Thus, for a configuration $\mathbf{r}\in (\mathbb{R}^2)^N$ and a  complex number $\textsf{z}\in \mathbb{C}$, $\textsf{z}\mathbf{r}$ is defined by
\begin{center}
$\textsf{z}\mathbf{r}:=(\textsf{z}\mathbf{r}_1, \textsf{z}\mathbf{r}_2, \cdots, \textsf{z}\mathbf{r}_N)$,
\end{center}
 and $\mathbf{i}\mathbf{r}$ is just $\mathbf{r}$ rotated anticlockwise by an angle $\frac{\pi}{2}$. Here and below,  please refer to \cite{yu2019problem} for more detail.
In this paper, unless otherwise specified, every vector is considered as a column vector.

For two  configurations $\mathbf{r} = (\mathbf{r}_1,  \cdots, \mathbf{r}_N)$ and $\mathbf{s} = (\mathbf{s}_1,  \cdots, \mathbf{s}_N)$ in $(\mathbb{R}^2)^N$, their \emph{mass scalar product}  is defined as:
\begin{displaymath}
\langle\mathbf{r},\mathbf{s}\rangle = \sum_{j=1}^{N} {{m_j (\mathbf{r}_j,\mathbf{s}_j)}}= \mathbf{r}^\top\mathfrak{M}\mathbf{s},
\end{displaymath}
where   $(\cdot,\cdot)$ denotes the standard scalar product in $\mathbb{R}^2$, $\mathfrak{M}$ is the diagonal matrix 
 \begin{center}
 $diag(m_1,m_1,m_2,m_2, \cdots, m_N, m_N)$,
\end{center}
and ``$\top$" denotes transposition of matrix. We denote $\|\cdot\|$ the Euclidean norm associated to the mass scalar product, that is
\begin{displaymath}
\|\mathbf{r}\|=\sqrt{\langle\mathbf{r},\mathbf{r}\rangle}.
\end{displaymath}
Then the cartesian space $(\mathbb{R}^2)^N $ is a new Euclidean space.

Recall that, the moment
of inertia, the kinetic energy, the opposite of the potential energy (force function), the total energy, the angular momentum, and the Lagrangian  function    are respectively defined as
\begin{displaymath}
I(\mathbf{r}) = \sum_{j=1}^{N} {{m_j |{\mathbf{r}}_j-{\mathbf{r}}_c|^2}}.
\end{displaymath}
\begin{displaymath}
\mathcal{K} (\dot{\mathbf{r}})= \sum_{j=1}^{N} {\frac{1}{2}{m_j |\dot{\mathbf{r}}_j|^2}},
\end{displaymath}
\begin{displaymath}
\mathcal{U}(\mathbf{r}) = \sum_{1\leq k<j\leq N} {\frac{m_k m_j }{|\mathbf{r}_k - \mathbf{r}_j|}},
\end{displaymath}
\begin{displaymath}
\mathcal{H }(\mathbf{r},\dot{\mathbf{r}})= \mathcal{K}(\dot{\mathbf{r}})- \mathcal{U}(\mathbf{r}),
\end{displaymath}
\begin{displaymath}
\mathcal{J}(\mathbf{r}) = \sum_{j=1}^{N} {m_j {\mathbf{r}}_j \times {\dot{\mathbf{r}}}_j},
\end{displaymath}
\begin{displaymath}
\mathcal{L}(\mathbf{r},\dot{\mathbf{r}}) = \mathcal{L} = \mathcal{K} + \mathcal{U} = \sum_j \frac{1}{2} m_j |\dot{\mathbf{r}}_j|^2  + \sum_{k<j}{\frac{m_k m_j}{|\mathbf{r}_k - \mathbf{r}_j|}},
\end{displaymath}
where $|\cdot|$  denotes the Euclidean norm in ${\mathbb{R}}^2$, $\times$ denotes the standard cross product in $\mathbb{R}^2$, and $\mathbf{r}_c = \frac{\sum_{k = 1}^{N} {m_k \vec{r}_k}}{\sum_{k = 1}^{N} {m_k }}$ is the center of mass.

Without loss of generality, one  assumes that the center of mass is fixed at the origin.
Let $\mathcal{X}$ denote the space of configurations whose center of mass is at the origin; that is, $\mathcal{X} = \{ \mathbf{r} = (\mathbf{r}_1,\cdots, \mathbf{r}_N)\in (\mathbb{R}^2)^N: \sum_{k = 1}^{N} {m_k \mathbf{r}_k} = 0  \}$, or,
\begin{displaymath}
\mathcal{X} = \{ \mathbf{r} \in (\mathbb{R}^2)^N: \langle\mathbf{r},\mathcal{E}_{2N-1}\rangle = 0, \langle\mathbf{r},\mathcal{E}_{2N}\rangle = 0   \},
\end{displaymath}
where
\begin{displaymath}
\mathcal{E}_{2N-1} = (1,0,\cdots,1,0)^\top, \mathcal{E}_{2N} = (0,1,\cdots,0,1)^\top.
\end{displaymath}
Thus \begin{equation*}
  (\mathbb{R}^2)^N =span\{\mathcal{E}_{2N-1},\mathcal{E}_{2N}\}\oplus \mathcal{X}.
\end{equation*}
Note that, for a configuration $\mathbf{r} \in \mathcal{X}$, we have
\begin{displaymath}
\|\mathbf{r}\|=\sqrt{I(\mathbf{r})}.
\end{displaymath}
Let $\Delta $ be the collision set in $(\mathbb{R}^2)^N $, that is, $\Delta = \{ \mathbf{r} \in (\mathbb{R}^2)^N :  \mathbf{r}_j = \mathbf{r}_k~ \text{for some}~ j<k\}$. Then
the set $\mathcal{X} \backslash \Delta$ is the space of collision-free configurations.\\

Let us recall the important concept of the central configuration:
\begin{definition}
A configuration $\mathbf{r} \in \mathcal{X} \backslash \Delta$ is called a central configuration if there exists a constant $\lambda\in {\mathbb{R}}$ such that
\begin{equation}\label{centralconfiguration}
\sum_{j=1,j\neq k}^N \frac{m_jm_k}{|\mathbf{r}_j-\mathbf{r}_k|^3}(\mathbf{r}_j-\mathbf{r}_k)=-\lambda m_k\mathbf{r}_k,1\leq k\leq N.
\end{equation}
\end{definition}
The value of $\lambda$ in (\ref{centralconfiguration}) is uniquely determined by
\begin{equation}\label{lambda}
\lambda=\frac{\mathcal{U}(\mathbf{r})}{I(\mathbf{r})}.
\end{equation}

It is well known that a central configuration  is just a critical point of the \emph{normalized potential}  $\widetilde{\mathcal{U}}:=I^{\frac{1}{2}}\mathcal{U}$. Moreover, a central configuration $\mathbf{r}_0$ will be called \emph{nondegenerate}, if  the kernel of $D^2\widetilde{\mathcal{U}}(\mathbf{r}_0)$, the Hessian  of $\widetilde{\mathcal{U}}$ evaluated at $\mathbf{r}_0$, is exactly the plane $\mathcal{P}_{\mathbf{r}_0}$, where \begin{center}
 $\mathcal{P}_{\mathbf{r}_0}:=\{\textsf{z}\mathbf{r}_0|\textsf{z}\in \mathbb{C}\}=span\{\mathbf{r}_0,\mathbf{i}\mathbf{r}_0\}$.
\end{center}

Let $\mathcal{P}^{\bot}_{\mathbf{r}_0}$ be the orthogonal complement of $\mathcal{P}_{\mathbf{r}_0}$ in $\mathcal{X}$, that is,\begin{equation*}
    \mathcal{X}=\mathcal{P}_{\mathbf{r}_0}\oplus\mathcal{P}^{\bot}_{\mathbf{r}_0},
\end{equation*} then $\mathcal{P}^{\bot}_{\mathbf{r}_0}$, $\mathcal{P}_{\mathbf{r}_0}$ and $span\{\mathcal{E}_{2N-1},\mathcal{E}_{2N}\}$ are three invariant subspaces of  $D^2\widetilde{\mathcal{U}}(\mathbf{r}_0)$.
Note that $\mathcal{E}_{2N-3}:=\mathbf{r}_0$, $ \mathcal{E}_{2N-2}:=\mathbf{i}\mathbf{r}_0$,
 $\mathcal{E}_{2N-1}$ and $\mathcal{E}_{2N}$ are four eigenvectors of the Hessian $D^2\widetilde{\mathcal{U}}(\mathbf{r}_0)$.
 By considering eigenvectors of the Hessian $D^2\widetilde{\mathcal{U}}(\mathbf{r}_0)$, it follows that there is an orthogonal basis\begin{center}
	$\{\mathcal{E}_1, \mathcal{E}_2, \cdots, \mathcal{E}_{2N-5},\mathcal{E}_{2N-4}, \mathcal{E}_{2N-3},\mathcal{E}_{2N-2}\}$
\end{center}of $\mathcal{X}$ such that
 \begin{center}
	$\mathcal{P}^{\bot}_{\mathbf{r}_0}=span\{\mathcal{E}_1, \mathcal{E}_2, \cdots, \mathcal{E}_{2N-5},\mathcal{E}_{2N-4}\}$,
\end{center}
 Assume
\begin{equation}\label{eigenorthogonalbasis}
    D^2\widetilde{\mathcal{U}}(\mathbf{r}_0)\mathcal{E}_j=\lambda_j\mathcal{E}_j,~~~~~~~~~~~~j=1,\cdots,2N,
\end{equation}
then $\lambda_{2N-3}=\lambda_{2N-2}=0$ and $\lambda_{2N-1}=\lambda_{2N}=I^{\frac{1}{2}} \lambda$.

 As a matter of fact, a straightforward computation shows that the Hessian $D^2\widetilde{\mathcal{U}}(\mathbf{r}_0)$ with
respect to the mass scalar product in $(\mathbb{R}^2)^N $
 is
\begin{displaymath}
I^{\frac{1}{2}} (\lambda \mathbb{I} + \mathfrak{M}^{-1} \mathfrak{B}) - 3 I^{-\frac{1}{2}} \lambda \mathbf{r}_0\mathbf{r}_0^\top \mathfrak{M},
\end{displaymath}
thus
\begin{equation}\label{eigenorthogonalbasis1}
    \left(I^{\frac{1}{2}} (\lambda \mathbb{I} + \mathfrak{M}^{-1} \mathfrak{B}) - 3 I^{-\frac{1}{2}} \lambda \mathbf{r}_0\mathbf{r}_0^\top \mathfrak{M}\right)\mathcal{E}_j=\lambda_j\mathcal{E}_j,~~~~~~~~~~~~j=1,\cdots,2N.
\end{equation}
Where $I=I(\mathbf{r}_0)=\|\mathbf{r}_0\|^2$, $\lambda=\frac{\mathcal{U}(\mathbf{r}_0)}{I(\mathbf{r}_0)}$, $\mathfrak{B}$ is the Hessian of $\mathcal{U}$ evaluated at $\mathbf{r}_0$ with respect to the standard scalar
product of $(\mathbb{R}^2)^N $ and can be viewed as an $N\times N$ array of $2 \times 2$ blocks:
\begin{center}
$\mathfrak{B} = \left(
       \begin{array}{ccc}
         B_{11} & \cdots & B_{1N} \\
         \vdots & \ddots & \vdots \\
         B_{N1} & \cdots & B_{NN}\\
       \end{array}
     \right).
$
\end{center}
The off-diagonal blocks are given by:
\begin{center}
$B_{jk} = \frac{m_j m_k}{r^3_{jk}}[\mathbb{I}-\frac{3(\mathbf{r}_k - \mathbf{r}_j)(\mathbf{r}_k - \mathbf{r}_j)^\top}{r^2_{jk}}],
$
\end{center}the diagonal blocks are given by:
\begin{displaymath}
B_{kk} = -\sum_{1\leq j\leq N, j\neq k} B_{jk},
\end{displaymath}
where $r_{jk}=|\mathbf{r}_k - \mathbf{r}_j|$. Note that, as a matter of notational convenience, the \emph{identity matrix} of any order will always be denoted by $\mathbb{I}$, and the order of $\mathbb{I}$ can be determined according to the context.

By \eqref{eigenorthogonalbasis1}, it follows that
\begin{equation}\label{Hessian111}
  (\lambda \mathbb{I} + \mathfrak{M}^{-1} \mathfrak{B} )  (\mathcal{E}_{1},\cdots, \mathcal{E}_{2N})
=(\mathcal{E}_{1},\cdots, \mathcal{E}_{2N}) diag(  \frac{\lambda_1}{\|\mathbf{r}_0\|}, \cdots,  \frac{\lambda_{2N-4}}{\|\mathbf{r}_0\|},3  \lambda, 0,\lambda, \lambda ).
\end{equation}
That is, the  orthogonal basis
	$\{\mathcal{E}_1,  \cdots,  \mathcal{E}_{2N}\}$
 and the corresponding eigenvalues $\lambda_1,  \cdots,  \lambda_{2N}$ can be obtained by calculating eigenvectors of the matrix $\lambda \mathbb{I} + \mathfrak{M}^{-1} \mathfrak{B}$.

Given a configuration $\mathbf{r}$, let $\hat{\mathbf{r}}:= \frac{\mathbf{r}}{\|\mathbf{r}\|}$  be the unit vector corresponding to $\mathbf{r}$ henceforth. In particular, the unit vector $\hat{\mathbf{r}}$ is called the  \emph{normalized configuration} of the configuration $\mathbf{r}$.

Then
\begin{displaymath}
\{\hat{\mathcal{E}}_1, \hat{\mathcal{E}}_2, \hat{\mathcal{E}}_3, \hat{\mathcal{E}}_4, \cdots, \hat{\mathcal{E}}_{2N}\}
\end{displaymath}
consisting of eigenvectors of $D^2\widetilde{\mathcal{U}}(\mathbf{r}_0)$ is a standard orthogonal basis of the space $(\mathbb{R}^2)^N $ with respect to the scalar product $\langle,\rangle$, that is,
\begin{displaymath}
(\hat{\mathcal{E}}_1, \hat{\mathcal{E}}_2, \hat{\mathcal{E}}_3, \hat{\mathcal{E}}_4, \cdots, \hat{\mathcal{E}}_{2N})^\top \mathfrak{M} (\hat{\mathcal{E}}_1, \hat{\mathcal{E}}_2, \hat{\mathcal{E}}_3, \hat{\mathcal{E}}_4, \cdots, \hat{\mathcal{E}}_{2N}) = \mathbb{I}.
\end{displaymath}

\subsection{Moving  Coordinates}
\indent\par
Let us recall the so-called moving   coordinates introduced in \cite{yu2019problem}.  The  moving   coordinates are appropriate to describe   the  orbits near a central configuration.
 Here we say a  configuration $\mathbf{r}\in \mathcal{X}$ is near the given central configuration $\mathbf{r}_0$, if $\mathbf{r}\notin\mathcal{P}^{\bot}_{\mathbf{r}_0}$.

For any configuration $\mathbf{r}\in \mathcal{X} \backslash  \mathcal{P}^{\bot}_{\mathbf{r}_0}$, it is easy to see that  there exists a unique point $e^{\mathbf{i}\theta(\mathbf{r})}\hat{\mathbf{r}}_0$ on the  circle $\mathbf{S}:= \{e^{\mathbf{i}\theta}\hat{\mathbf{r}}_0| \theta \in \mathbb{R}\}$ such that
\begin{equation}
\|e^{\mathbf{i}\theta(\mathbf{r})}\hat{\mathbf{r}}_0-\mathbf{r}\|= min_{\theta \in \mathbb{R}} \|e^{\mathbf{i}\theta}\hat{\mathbf{r}}_0 - \mathbf{r}\|,\nonumber
\end{equation}
 Indeed,  the angle $\theta(\mathbf{r})$, unique up to integer multiple of $2\pi$, is determined by
 $\langle e^{-\mathbf{i}\theta(\mathbf{r})}\hat{\mathbf{r}},\mathbf{i} \hat{\mathbf{r}}_0\rangle=0$
 and $\langle e^{-\mathbf{i}\theta(\mathbf{r})}\hat{\mathbf{r}}, \hat{\mathbf{r}}_0\rangle>0$.

By decomposing $ e^{-\mathbf{i}\theta(\mathbf{r})}\hat{\mathbf{r}}$ with respect to the basis $\{\hat{\mathcal{E}}_{1}, \cdots, \hat{\mathcal{E}}_{2N-2}\}$,
and denoting $z=(z_1, \cdots, z_{2N-4})^\top$ the coordinates with respect to the vectors $\hat{\mathcal{E}}_{1}, \cdots, \hat{\mathcal{E}}_{2N-4}$, it holds
\begin{equation}\label{movingcoordinates}
    \mathbf{r}=r \hat{\mathbf{r}}=re^{\mathbf{i}\theta(\mathbf{r})}({\sqrt{1 - |z|^2} \hat{\mathbf{r}}_0}+\sum_{k = 1}^{2N-4} {z_k \hat{\mathcal{E}}_k});
\end{equation}
where
 $|z|^2=z^\top z=\sum_{j = 1}^{2N-4} z^2_j$.

Then the total set of the variables $r, \theta, z$  are referred as the \emph{moving  coordinates} of $\mathbf{r}\in \mathcal{X} \backslash  \mathcal{P}^{\bot}_{\mathbf{r}_0}$.

Note that
\begin{displaymath}
min_{\theta \in \mathbb{R}} \|e^{\mathbf{i}\theta}\hat{\mathbf{r}}_0 - \mathbf{r}\| \geq 1,
\end{displaymath}
if  $\mathbf{r}\in   \mathcal{P}^{\bot}_{\mathbf{r}_0}$. Consequently, if
\begin{displaymath}
min_{\theta \in \mathbb{R}} \|e^{\mathbf{i}\theta}\hat{\mathbf{r}}_0 - \mathbf{r}\| < 1,
\end{displaymath}
then  $\mathbf{r}\in  \mathcal{X} \backslash \mathcal{P}^{\bot}_{\mathbf{r}_0}$. Therefore we have legitimate rights to use the coordinates $ (r, \theta, z)$ in a neighbourhood of $\mathbf{r}_0$.

\subsection{Equations of Motion in the Moving  Coordinates}
\indent\par

Now let us write the equations of motion in the moving coordinates.

It is well known that the equations (\ref{eq:Newton's equation1}) of motion are just the Euler-Lagrange equations
of  Lagrangian
system with the configuration manifold $\mathcal{X}$ and the Lagrangian function $\mathcal{L}(\mathbf{r},\dot{\mathbf{r}})$. Then this yields a quick method
for writing equations of motion in various coordinate systems. In fact, to write the equations of motion in a new coordinate system, it is sufficient to
express the Lagrangian function in the new coordinates. Please refer to the appendix for more detail.

By \begin{equation*}
    \mathbf{r}=r e^{\mathbf{i}\theta}(z_0 \hat{\mathbf{r}}_0+\sum_{k = 1}^{2N-4} {z_k \hat{\mathcal{E}}_k}),
\end{equation*}where $z_0=\sqrt{1 - |z|^2}$, it follows that the kinetic energy and the force function can be written as
\begin{equation}\label{kineticenergy}
\mathcal{K}(\mathbf{r}) =\frac{1}{2}\langle\dot{\mathbf{r}},\dot{\mathbf{r}}\rangle=\frac{\dot{r}^2}{2}+ \frac{r^2}{2} (\dot{z}^2_0+|\dot{z}|^2 + 2\dot{\theta}\dot{z}^\top Q {z}+ \dot{\theta}^2),\nonumber
\end{equation}
\begin{equation}\label{force function}
\mathcal{U}(\mathbf{r}) = \frac{\mathcal{U}(z_0 \hat{\mathbf{r}}_0+\sum_{j = 1}^{2N-4} {z_j \hat{\mathcal{E}}_j})}{r},\nonumber
\end{equation}
where
the square matrix
\begin{equation*}
    Q:=(q_{jk})_{(2N-4)\times (2N-4)}
\end{equation*}
 is an anti-symmetric orthogonal matrix such that $q_{jk}=\langle\hat{\mathcal{E}}_j,\mathbf{i}\hat{\mathcal{E}}_k\rangle$. Note that
 \begin{equation*}
\mathcal{U}(z_0 \hat{\mathbf{r}}_0+\sum_{j = 1}^{2N-4} {z_j \hat{\mathcal{E}}_j})=\widetilde{\mathcal{U}}(z_0\hat{\mathbf{r}}_0+\sum_{j = 1}^{2N-4} {z_j \hat{\mathcal{E}}_j})
 \end{equation*}
 only contains the variables $z_j$ ($j=1, \cdots, 2N-4$), we will simply
write it as  $U(z)$ henceforth. 
Therefore, \begin{equation}\label{force function1}
\mathcal{U}(\mathbf{r}) = \frac{U(z)}{r}.\nonumber
\end{equation}

As a result,  the Lagrangian  function $\mathcal{L}$ is
\begin{equation}\label{Lagrangianfunction1}
    \mathcal{L}(z,r,\dot{z},\dot{r},\dot{\theta})=\frac{\dot{r}^2}{2}+ \frac{r^2}{2} (\dot{z}^2_0+|\dot{z}|^2 + 2\dot{\theta}\dot{z}^\top Q {z}+ \dot{\theta}^2)+\frac{U(z)}{r}.
\end{equation}
The equations of motion in the moving coordinates are just the Euler-Lagrange equations
of (\ref{Lagrangianfunction1}):
\begin{equation*}\label{Euler-Lagrangeequations0}
\left\{
             \begin{array}{l}
           \frac{d}{dt}\frac{\partial \mathcal{L}}{\partial \dot{z}}-\frac{\partial \mathcal{L}}{\partial z}=0,    \\
             \frac{d}{dt}\frac{\partial \mathcal{L}}{\partial \dot{r}}-\frac{\partial \mathcal{L}}{\partial r}=0,   \\
            \frac{d}{dt}\frac{\partial \mathcal{L}}{\partial \dot{\theta}}-\frac{\partial \mathcal{L}}{\partial \theta}=0;
             \end{array}
\right.
\end{equation*}
or\begin{equation}\label{equations of motion}
\left\{
             \begin{array}{l}
             {r}^2 [ \ddot{z} + \ddot{\theta} Qz+ 2 \dot{\theta} Q \dot{z}]
             +{r}^2 [\frac{ z^\top\ddot{z}+ |\dot{z}|^2}{1 -  |z|^2} + \frac{ ({z}^\top \dot{z})^2}{(1 -  |z|^2)^2}]z
 + 2r\dot{r}[\frac{ {z}^\top \dot{z}}{1 -  |z|^2}z + \dot{z} -  \dot{\theta}Q z] - \frac{1}{r}\frac{\partial U(z)}{\partial z}=0,   \\
             \ddot{r}- {r} [\frac{ ({z}^\top \dot{z})^2}{1 -  |z|^2}+|\dot{z}|^2 + 2\dot{\theta}\dot{z}^\top Q {z}+ \dot{\theta}^2]+ \frac{ U(z)}{r^2}=0,   \\
             2r\dot{r} (\dot{z}^\top Q {z}+ \dot{\theta}) + r^2 (\ddot{z}^\top Q {z}+ \ddot{\theta})=0.
             \end{array}
\right.
\end{equation}

It is noteworthy that the variable $\theta$ is not involved in the  Lagrangian
$\mathcal{L}$,
this is a main reason of introducing the moving  coordinates.
That is, the variable $\theta$ is an ignorable coordinate.
As a result, $\frac{\partial \mathcal{L}}{\partial \dot{\theta}}$ is conserved. In fact,   a straightforward computation shows that
\begin{equation}\label{angular momentum}
\frac{\partial \mathcal{L}}{\partial \dot{\theta}}= \mathcal{J} = \langle{\mathbf{i}\mathbf{r}},{\dot{\mathbf{r}}}\rangle
= r^2(\dot{\theta}+\dot{z}^\top Q {z}).
\end{equation}

By \emph{Routh's method for eliminating
ignorable coordinates}, we introduce the function
\begin{equation}\label{reducedLagrangian}
\mathcal{L}_{\mathcal{J}}(z,r,\dot{z},\dot{r})=\mathcal{L}(z,r,\dot{z},\dot{r},\dot{\theta}) - {\mathcal{J}} \dot{\theta}|_{r,z,\dot{r},\dot{z},{\mathcal{J}}} \nonumber
\end{equation}
as the reduced Lagrangian function on the level set of $\mathcal{J}$, where $\mathcal{J}$ is certainly a constant and we represent
$\dot{\theta}$ as a function of $z,r,\dot{z},\dot{r}$ and $\mathcal{J}$ by using the equality (\ref{angular momentum}). A straightforward computation shows that
\begin{equation}\label{reducedLagrangian1}
\mathcal{L}_{\mathcal{J}}(z,r,\dot{z},\dot{r})= \frac{\dot{r}^2}{2}+ \frac{r^2}{2} \left[\dot{z}^2_0+|\dot{z}|^2 -\left(\dot{z}^\top Q {z}- \frac{\mathcal{J}}{r^2}\right)^2\right]+\frac{U(z)}{r}.
\end{equation}
The equations of motion on the level set of $\mathcal{J}$ are just the Euler-Lagrange equations
of (\ref{reducedLagrangian1}):
\begin{equation*}\label{Euler-Lagrangeequations1}
\left\{
             \begin{array}{l}
           \frac{d}{dt}\frac{\partial \mathcal{L}_{\mathcal{J}}}{\partial \dot{z}}-\frac{\partial \mathcal{L}_{\mathcal{J}}}{\partial z}=0,    \\
             \frac{d}{dt}\frac{\partial \mathcal{L}_{\mathcal{J}}}{\partial \dot{r}}-\frac{\partial \mathcal{L}_{\mathcal{J}}}{\partial r}=0;
             \end{array}
\right.
\end{equation*}
or
\begin{equation}\label{Euler-Lagrangeequations1detail}
\left\{
             \begin{array}{lr}
              {r}^2 [\ddot{z} - (\ddot{z}^\top Q {z})Qz- 2  (\dot{z}^\top Q {z})Q \dot{z}]
               +{r}^2 [\frac{ z^\top\ddot{z}+ |\dot{z}|^2}{1 -  |z|^2} + \frac{ ({z}^\top \dot{z})^2}{(1 -  |z|^2)^2}]z+ 2\mathcal{J}Q \dot{z}\\
 +2r\dot{r}[\frac{ {z}^\top \dot{z}}{1 -  |z|^2}z + \dot{z} -  (\dot{z}^\top Q {z})Q z] - \frac{1}{r}\frac{\partial U(z)}{\partial z_i}=0,  & \\
             \ddot{r}- {r} [\frac{ ({z}^\top \dot{z})^2}{1 -  |z|^2}+|\dot{z}|^2  - (\dot{z}^\top Q {z})^2+\frac{\mathcal{J}^2}{r^4}]+ \frac{ U(z)}{r^2}=0.  &
             \end{array}
\right.
\end{equation}

It is noteworthy that there is no degeneracy  according to translation symmetry and   rotation symmetry    of the $N$-body problem  in the system (\ref{Euler-Lagrangeequations1detail}).

\section{ Stability on Relative Equilibria}
\indent\par

For the central configuration  $\mathbf{r}_0$, let us consider a
 relative equilibrium $\rho  e^{\mathbf{i}\omega t}\mathbf{r}_0$ of the Newtonian $N$-body problem, where  $\rho$ and $\omega$ are two positive constants. Without
loss of generality, suppose that $\|\mathbf{r}_0\|=1$, i.e., $\mathbf{r}_0=\hat{\mathbf{r}}_0$.
 Then a straight forward computation shows that the angular momentum $\mathcal{J}$ of the relative equilibrium $\rho  e^{\mathbf{i}\omega t}\mathbf{r}_0$ is just $\omega \rho^2$ and $\lambda = \rho^3 \omega^2$.

By the moving coordinates $ r, \theta, z$, the relative equilibrium $\rho  e^{\mathbf{i}\omega t}\mathbf{r}_0$ is just a solution of the system (\ref{equations of motion}) such that
\begin{displaymath}
r=\rho, \theta=\omega t, z=0.
\end{displaymath}
In the system (\ref{Euler-Lagrangeequations1detail}), the relative equilibrium $\rho  e^{\mathbf{i}\omega t}\mathbf{r}_0$ is just corresponding to an equilibrium solution  such that
\begin{displaymath}
r=\rho,  z=0.
\end{displaymath}

Recall that, every elliptic  solution of the two-body problem is always instable in the sense of Lyapunov, but  is orbitally  stable.
Similar to the case of the two-body problem,
  we have the fact that the relative equilibrium $\rho  e^{\mathbf{i}\omega t}\mathbf{r}_0$ is always  instable in the sense of Lyapunov. Here we omit  the proof of the fact and only remind that the  proof   attributes to the discussion of the two-body problem.
Therefore, it is  natural to consider the orbital  stability of  the relative equilibrium $\rho  e^{\mathbf{i}\omega t}\mathbf{r}_0$. 

\subsection{ Orbital Stability}
\indent\par
Let $\mathbf{S}_{\rho}$, a circle, denote the orbit of the relative equilibrium $\rho  e^{\mathbf{i}\omega t}\mathbf{r}_0$ in the phase space of the $N$-body problem.
By introducing new variables:
\begin{displaymath}
Z=\dot{z},  ~~~~~~~~~~\Upsilon = {\dot{r}}, ~~~~~~~~~~\Theta = \dot{\theta},
\end{displaymath}
 the system (\ref{equations of motion}) becomes:
\begin{equation}\label{equations of motion0}
\left\{
             \begin{array}{lr}
             \dot{z} = Z, ~~~~~~~ & \\
             \dot{Z} =- \dot{\Theta} Qz- 2 \Theta Q Z
             - [\frac{ z^\top\dot{Z}+ |Z|^2}{1 -  |z|^2} + \frac{ ({z}^\top Z)^2}{(1 -  |z|^2)^2}]z 
 - \frac{2\Upsilon }{r}[\frac{ {z}^\top Z}{1 -  |z|^2}z + Z -  \Theta Q z]
 + \frac{1}{r^3}\frac{\partial U(z)}{\partial z},&\\
             {\dot{r}}=\Upsilon, &\\
             \dot{\Upsilon}= r [\frac{ ({z}^\top Z)^2}{1 -  |z|^2}+|Z|^2 + 2\Theta Z^\top Q {z}+ \Theta^2]- \frac{ U(z)}{r^2}, &\\
             {\dot{\theta}}=\Theta, &\\
              \dot{\Theta}=-\dot{Z}^\top Q {z}- \frac{2\Upsilon}{r}(Z^\top Q {z}+ \Theta  ) , &
             \end{array}
\right.
\end{equation}
and the system (\ref{Euler-Lagrangeequations1detail}) becomes:
\begin{equation}\label{equations of motion00}
\left\{
             \begin{array}{lr}
              \dot{z} = Z, ~~~~~~~ & \\
             \dot{Z} = (\dot{Z}^\top Q {z})Qz+ 2  (Z^\top Q {z})Q Z
               - [\frac{ z^\top\dot{Z}+ |Z|^2}{1 -  |z|^2} + \frac{ ({z}^\top Z)^2}{(1 -  |z|^2)^2}]z- \frac{2\mathcal{J}}{{r}^2}Q Z\\
 -\frac{2\Upsilon}{{r}}[\frac{ {z}^\top Z}{1 -  |z|^2}z + Z -  (Z^\top Q {z})Q z] + \frac{1}{{r}^3}\frac{\partial U(z)}{\partial z_i},  & \\[8pt]
             {\dot{r}}=\Upsilon, &\\
             \dot{\Upsilon}= {r} [\frac{ ({z}^\top Z)^2}{1 -  |z|^2}+|Z|^2  - (Z^\top Q {z})^2+\frac{\mathcal{J}^2}{r^4}]- \frac{ U(z)}{r^2}.  &
             \end{array}
\right.
\end{equation}
Then the relative equilibrium $\rho  e^{\mathbf{i}\omega t}\mathbf{r}_0$ is just corresponding to a solution of the system (\ref{equations of motion0}) such that
\begin{displaymath}
z=0, Z=0, r=\rho, \Upsilon=0, \theta=\omega t, \Theta=\omega;
\end{displaymath}
 and just corresponding to an equilibrium point of the system (\ref{equations of motion00})   such that
\begin{displaymath}
z=0, Z=0, r=\rho, \Upsilon=0.
\end{displaymath}

By the concept of classical orbital  stability, we define orbital  stability of the relative equilibrium $\rho  e^{\mathbf{i}\omega t}\mathbf{r}_0$:
\begin{definition}
The relative equilibrium $\rho  e^{\mathbf{i}\omega t}\mathbf{r}_0$ is orbitally stable,
if  for any $\epsilon>0$, there is a neighbourhood $\mathcal{N}$ of the orbit $\mathbf{S}_{\rho}$,
such that,  the distance
from any orbit, whose initial value belongs to $\mathcal{N}$, to $\mathbf{S}_{\rho}$ is less than $\epsilon$.
\end{definition}
We claim that
\begin{theorem}\label{Orbitalstability of equations of motion}
The relative equilibrium $\rho A(\omega t)\mathcal{E}_3$ is orbitally stable  if and only if the equilibrium point $z=0, Z=0, r=\rho, \Upsilon=0$ of the system (\ref{equations of motion00}) is stable in the sense of Lyapunov.
\end{theorem}
{\bf Proof of Theorem \ref{Orbitalstability of equations of motion}:}

Recall the definition of  the  moving coordinates, it is clear that the quantity 
\begin{center}
$|z|^2 + |Z|^2 + (r-\rho)^2+\Upsilon^2+(\Theta-\omega)^2$
\end{center}
 gives a measure of the distance
from a phase point $(z,Z,r,\Upsilon,\theta,\Theta)$ to the orbit $\mathbf{S}_{\rho}$. 

If we replace the variable $\Theta=\dot{\theta}$ by the variable $\vartheta=\frac{\partial \mathcal{L}}{\partial \dot{\theta}}$, then the variables $z,Z,r,\Upsilon,\theta,\vartheta$ are new  coordinates of phase points. Note that the variable $\vartheta$ is  just $\mathcal{J}$, and by (\ref{angular momentum}) (\ref{Euler-Lagrangeequations1detail}),  the system (\ref{equations of motion0}) becomes
\begin{equation}\label{equations of motion0n}
\left\{
             \begin{array}{lr}
              \dot{z} = Z, ~~~~~~~ & \\
             \dot{Z} = (\dot{Z}^\top Q {z})Qz+ 2  (Z^\top Q {z})Q Z
               - [\frac{ z^\top\dot{Z}+ |Z|^2}{1 -  |z|^2} + \frac{ ({z}^\top Z)^2}{(1 -  |z|^2)^2}]z- \frac{2\mathcal{J}}{{r}^2}Q Z\\
 -\frac{2\Upsilon}{{r}}[\frac{ {z}^\top Z}{1 -  |z|^2}z + Z -  (Z^\top Q {z})Q z] + \frac{1}{{r}^3}\frac{\partial U(z)}{\partial z_i},  & \\[6pt]
             {\dot{r}}=\Upsilon, &\\
             \dot{\Upsilon}= {r} [\frac{ ({z}^\top Z)^2}{1 -  |z|^2}+|Z|^2  - (Z^\top Q {z})^2+\frac{\mathcal{J}^2}{r^4}]- \frac{ U(z)}{r^2},  &\\
             \dot{\theta}=\frac{ \mathcal{J}}{r^2}-Z^\top Q {z},  &\\
              \dot{\vartheta}=0.  &
             \end{array}
\right.
\end{equation}
It is obvious  that the relative equilibrium $\rho  e^{\mathbf{i}\omega t}\mathbf{r}_0$ is just corresponding to a solution of the system (\ref{equations of motion0n}) such that
\begin{displaymath}
z=0, Z=0, r=\rho, \Upsilon=0, \theta=\omega t, \vartheta=\omega \rho^2;
\end{displaymath}
and the quantity
\begin{center}
$|z|^2 + |Z|^2 + (r-\rho)^2+\Upsilon^2+(\vartheta-\omega \rho^2)^2$
\end{center}
 gives a measure of the distance
from a phase point $(z,Z,r,\Upsilon,\theta,\vartheta)$ to the orbit $\mathbf{S}_{\rho}$. 

Note that, by $\dot{\vartheta}=0$, $(\vartheta-\omega \rho^2)^2$ is always small provided that its initial value is small.
As a result, to prove orbital  stability of the relative equilibrium $\rho  e^{\mathbf{i}\omega t}\mathbf{r}_0$, it  is equivalent to show that in the system (\ref{equations of motion0n})
the quantity
\begin{center}
$|z|^2 + |Z|^2 + (r-\rho)^2+\Upsilon^2$
\end{center}
is always small if its initial value is small. It is easy to see that it  is equivalent to show that the equilibrium point 
 \begin{center}
 $z=0, Z=0, r=\rho, \Upsilon=0$
\end{center}
of the system (\ref{equations of motion00}) is stable in the sense of Lyapunov.


$~~~~~~~~~~~~~~~~~~~~~~~~~~~~~~~~~~~~~~~~~~~~~~~~~~~~~~~~~~~~~~~~~~~~~~~~~~~~~~~~~~~~~~~~~~~~~~~~~~~~~~~~~~~~~~~~~~~~~~~~~~~~~~~~~~~~~~~~~~~~~~~~~~~~~~~~~~~~~~~~~~~\Box$\\

As a matter of convenience, the equilibrium point $z=0, Z=0, r=\rho, \Upsilon=0$ in the system (\ref{equations of motion00}) will be translated to the origin by substituting $r+\rho$ for $r$. Then the problem is now reduced to investigate 
stability of the origin of the following system
\begin{equation}\label{equations of motion1}
\left\{
             \begin{array}{lr}
              \dot{z} = Z, ~~~~~~~ & \\
             \dot{Z} = (\dot{Z}^\top Q {z})Qz+ 2  (Z^\top Q {z})Q Z
               - [\frac{ z^\top\dot{Z}+ |Z|^2}{1 -  |z|^2} + \frac{ ({z}^\top Z)^2}{(1 -  |z|^2)^2}]z- \frac{2\mathcal{J}}{{(r+\rho)}^2}Q Z\\
 -\frac{2\Upsilon}{r+\rho}[\frac{ {z}^\top Z}{1 -  |z|^2}z + Z -  (Z^\top Q {z})Q z] + \frac{1}{{(r+\rho)}^3}\frac{\partial U(z)}{\partial z_i},  & \\[8pt]
             {\dot{r}}=\Upsilon, &\\
             \dot{\Upsilon}= {(r+\rho)} [\frac{ ({z}^\top Z)^2}{1 -  |z|^2}+|Z|^2  - (Z^\top Q {z})^2+\frac{\mathcal{J}^2}{(r+\rho)^4}]- \frac{ U(z)}{(r+\rho)^2}.  &
             \end{array}
\right.
\end{equation}

\subsection{Linear   Stability}
\indent\par
First, let us investigate   linear stability of the origin of the  system (\ref{equations of motion1}). Although our method has some differences from the classical works in \cite[etc]{moeckel1995linear,roberts1999spectral}, there is no new result for linear stability. So we only consider some special cases to reveal the good of the moving coordinates.

Linearizing the  system \eqref{equations of motion1} at the  origin yields the following linearized system
\begin{equation*}
\left\{
             \begin{array}{lr}
             \dot{z} = Z, &  \\
             \dot{Z} = \frac{ \Lambda z }{\rho^3} - 2 \omega Q Z,  & \\
             {\dot{r}}=\Upsilon, &\\
             \dot{\Upsilon}= -\omega^2r; &
             \end{array}
\right.
\end{equation*}
or
\begin{equation}\label{linearized equations1}
\left(
  \begin{array}{c}
    \dot{z} \\
    \dot{Z} \\
    \dot{r} \\
    \dot{\Upsilon} \\
  \end{array}
\right)=\left(
          \begin{array}{cccc}
            0  & \mathbb{I} &   &      \\
            \frac{ \Lambda}{\rho^3} & -2\omega Q &   &     \\
              &   & 0  & 1  \\
              &   & -\omega^2 & 0 \\
          \end{array}
        \right)\left(
\begin{array}{c}
    {z} \\
    {Z} \\
    {r} \\
    {\Upsilon} \\
  \end{array}
\right),
\end{equation}
where
$\Lambda = diag( \lambda_1, \cdots,  \lambda_{2N-4})$. In the calculation, please note that $\mathcal{J}=\omega \rho^2$, $\lambda = \rho^3 \omega^2$,
and, as in \cite{yu2019problem}, we can expand $U(x)$ as\begin{equation}
U(z)= \lambda + \frac{1}{2}\sum_{k = 1}^{2N-4} \lambda_k z^2_k   +\cdots,\nonumber
\end{equation}where ``$\cdots$" denotes the terms of degree higher than 2.

Following  Moeckel's approach in \cite{moeckel1995linear}, we define linear stability and spectral stability of the relative equilibrium $\rho  e^{\mathbf{i}\omega t}\mathbf{r}_0$:
\begin{definition}
The relative equilibrium $\rho  e^{\mathbf{i}\omega t}\mathbf{r}_0$ is spectrally stable, 
if all the eigenvalues of the matrix in the linearized system (\ref{linearized equations1}) are either zero or purely imaginary. The relative equilibrium $\rho  e^{\mathbf{i}\omega t}\mathbf{r}_0$ is called linearly stable, if it is  spectrally stable and the matrix in  (\ref{linearized equations1}) is further diagonalizable.
\end{definition}

If a solution is   spectrally instable, it follows from the well known Lyapunov's theorem of stability that the solution is Lyapunov  instable. However, if a solution is  spectrally stable but  linearly instable, it is possible that the solution is Lyapunov stable; similarly, if a solution is  linearly stable, it is possible that the solution is Lyapunov instable.

It is easy to see that a Jordan canonical form of the submatrix $\begin{pmatrix}
 0  & 1  \\
  -\omega^2 & 0  \\
          \end{pmatrix}$ is $
          \begin{pmatrix}
              \mathbf{i}\omega  & 0 \\
               0 & -\mathbf{i}\omega \\
          \end{pmatrix}$. It follows that we have
\begin{theorem}\label{stability of linearized equations1}
The relative equilibrium $\rho  e^{\mathbf{i}\omega t}\mathbf{r}_0$ is spectrally  stable if and only if all the eigenvalues of the matrix $\begin{pmatrix}
            0  & \mathbb{I} \\
            \frac{ \Lambda}{\rho^3} & -2\omega Q \\
          \end{pmatrix}$ are either zero or purely imaginary; and the relative equilibrium $\rho  e^{\mathbf{i}\omega t}\mathbf{r}_0$ is  linearly stable if and only if the matrix $\begin{pmatrix}
            0  & \mathbb{I} \\
            \frac{ \Lambda}{\rho^3} & -2\omega Q \\
          \end{pmatrix}$ is diagonalizable and all of its eigenvalues  are either zero or purely imaginary. 
\end{theorem}

Without loss of generality, assume $\rho=1$ from now on. Then $\mathcal{J}=\omega $ and $\lambda = \omega^2$.\\

\textbf{A. General case.}  For linear stability and spectral stability, it is necessary that the roots of the characteristic polynomial
\begin{displaymath}
\left|
          \begin{array}{cc}
            x \mathbb{I}  & -\mathbb{I}          \\
            -\Lambda & x \mathbb{I}-2\omega Q     \\
          \end{array}
        \right|=\left|x^2 \mathbb{I}-2\omega x Q -\Lambda
        \right|
\end{displaymath}
are all on the imaginary axis.
Due to $Q^\top=-Q$, the above characteristic polynomial is an even function of $x$. Let
\begin{displaymath}
f(x^2)=x^{2n} + c_{n-1}x^{2n-2}+\cdots+c_2 x^{4}+c_1 x^{2}+c_0
\end{displaymath}
be the above characteristic polynomial, where 
\begin{center}
$n:=2N-4$.
\end{center}
 Then it is necessary that all the roots of $f$ are either negative  or zero. It follows that
\begin{displaymath}
 c_{j}\geq 0 ~~~~~~~~~for~j=0,1,\cdots, n-1.
\end{displaymath}
Furthermore, if $c_{k}= 0$, then $c_{j}= 0$ for $0\leq j \leq k-1$.

A straightforward computation shows that:
\begin{displaymath}
\begin{array}{c}
 c_{n-1}  =2n\lambda- \sum^{n}_{j=1}\lambda_j,\\
c_0 =\Pi^{n}_{j=1}\lambda_j.
\end{array}
\end{displaymath}
It follows from (\ref{Hessian111}) that
\begin{displaymath}
\sum^{n}_{j=1}\lambda_j= (2N-5)\lambda  + \sum^{2N}_{k=1}\sum_{j\neq k} \frac{m_j}{r_{jk}^3}=(2N-5)\lambda  + \sum_{1\leq j< k\leq2N} \frac{m_j+m_k}{r_{jk}^3}.
\end{displaymath}
Thus
\begin{displaymath}
(2N-3)\lambda > \sum_{1\leq j< k\leq2N} \frac{m_j+m_k}{r_{jk}^3}.
\end{displaymath}
From the inequality above,   Roberts \cite{roberts1999spectral} proved that any relative equilibrium of $N$ equal masses is  spectrally
instable for $N \geq 24 306$.

\textbf{B. Collinear case.} When the central configuration $\mathbf{r}_0$ is collinear, the case becomes simpler. 

Suppose $\mathbf{r}_0 = (x_1, 0, x_3, 0, \cdots, x_{2N-1}, 0)^\top\in (\mathbb{R}\times 0)^{N} \subset \mathbb{R}^{2N}$, then the matrix $B_{jk} = \frac{m_j m_k}{r^3_{jk}}D$, where $D=\left(
                                    \begin{array}{cc}
                                      -2 & 0 \\
                                      0 & 1 \\
                                    \end{array}
                                  \right)$. So $\mathfrak{M}^{-1} \mathfrak{B}$ becomes:
\begin{center}
$ \left(
       \begin{array}{cccc}
         A_{11} &\frac{ m_2}{r^3_{12}} D &\cdots & \frac{ m_N}{r^3_{1N}} D \\
         \frac{ m_1}{r^3_{12}} D &A_{22} &\cdots & \frac{ m_N}{r^3_{2N}} D \\
         \vdots &\vdots & \ddots & \vdots \\
         \frac{ m_1}{r^3_{1N}} D & \frac{ m_2}{r^3_{2N}} D& \cdots & A_{NN}\\
       \end{array}
     \right),
$
\end{center}
where  the diagonal blocks are given by:
\begin{displaymath}
A_{kk} = -\sum_{1\leq j\leq N, j\neq k} \frac{ m_j}{r^3_{jk}} D.
\end{displaymath}

It follows that, if
\begin{displaymath}
(\lambda \mathbb{I}+\mathfrak{M}^{-1} \mathfrak{B})\mathcal{E}= \lambda_{\ast} \mathcal{E}  ~~~~~~~~~~~for ~\mathcal{E}\in (\mathbb{R}\times 0)^{N},
\end{displaymath}
then
\begin{displaymath}
(\lambda \mathbb{I}+\mathfrak{M}^{-1} \mathfrak{B})\mathbf{i}\mathcal{E}= \frac{3\lambda -\lambda_{\ast}}{2} \mathbf{i}\mathcal{E}.
\end{displaymath}
That is to say, if a vector $\mathcal{E}\in (\mathbb{R}\times 0)^{N}$ is an eigenvector of the Hessian $D^2\widetilde{\mathcal{U}}(\mathbf{r}_0)$ with eigenvalue $\|\mathbf{r}_0\|\lambda_{\ast}$, then $\mathbf{i}\mathcal{E} \in (0 \times \mathbb{R})^{N} \subset \mathbb{R}^{2N}$ is  an eigenvector of the Hessian $D^2\widetilde{\mathcal{U}}(\mathbf{r}_0)$ with eigenvalue $\frac{3\|\mathbf{r}_0\|\lambda -\|\mathbf{r}_0\|\lambda_{\ast}}{2}$.
Therefore, we can assume that the family $\{\mathcal{E}_1, \mathcal{E}_2,  \cdots, \mathcal{E}_{2N-5}, \mathcal{E}_{2N-4}\}$ satisfy
\begin{equation*}
   \mathcal{E}_{2k}=\mathbf{i}\mathcal{E}_{2k-1} 
\end{equation*}for $k=1,\cdots,N-2$.
Then $Q$ becomes block diagonal with  block $J=\left(
           \begin{array}{cc}
             0 & -1 \\
             1 & 0 \\
           \end{array}
         \right)$:
         \begin{equation*}
            Q = \left(
       \begin{array}{ccc}

         J &   &  \\
           & \ddots &   \\
          &   & J\\
       \end{array}
     \right).
         \end{equation*}
It follows that the characteristic polynomial $f(x^2)$ becomes
\begin{displaymath}
f(x^2)=\prod_{k=1}^{N-2}\left(x^{4} + (4\lambda-\lambda_{2k-1}-\lambda_{2k}) x^{2}+\lambda_{2k-1}\lambda_{2k}\right).
\end{displaymath}
According to a well known result  due to
Conley \cite{pacella1987central}, it holds  $\lambda_{2k-1}\lambda_{2k} <0$, thus the central configuration $\mathbf{r}_0$ is  spectrally instable. This result has been proved by
Moeckel \cite{moeckel1995linear}.
Furthermore, we have the following result of collinear central configurations:
\begin{equation}\label{lineareigenvalues}
\lambda_{2k}=\frac{3\lambda -\lambda_{2k-1}}{2}, ~~~~~~~\lambda_{2k-1} >3\lambda.
\end{equation}

So we  consider  stability of  relative equilibria only for noncollinear central configurations in the following.

\textbf{C. The three-body case.} When $N=3$, the problem is especially simple, because of $Q=\left(
                                                           \begin{array}{cc}
                                                             0 & 1 \\
                                                             -1 & 0 \\
                                                           \end{array}
                                                         \right)
$ or $\left(
        \begin{array}{cc}
          0 & -1 \\
          1 & 0 \\
        \end{array}
      \right)
$. Then all the roots of the characteristic polynomial
\begin{displaymath}
f(x^2)=x^{4} + (4\lambda-\lambda_1-\lambda_2) x^{2}+\lambda_1\lambda_2
\end{displaymath}
 are  on  the imaginary axis if and only if
\begin{equation}\label{threebodylinearstability}
\begin{array}{c}
 4\lambda\geq \lambda_1+\lambda_2,\\
 \lambda_1\lambda_2\geq 0, \\
 (4\lambda-\lambda_1-\lambda_2)^2-4\lambda_1\lambda_2\geq0.
\end{array}
\end{equation}

Set $\beta=m_1 m_2+ m_3 m_2+ m_1 m_3$. Since $\mathbf{r}_0$ is a equilateral triangle such that $\|\mathbf{r}_0\|=1$ and whose center of masses is at the origin,
without loss of generality, suppose
\begin{displaymath}
 m_1+m_2+m_3=1,
\end{displaymath}and
\begin{equation}\label{Lagrangetriangularpoint}
     \mathbf{r}_0= (-\frac{\sqrt{3} m_3}{2 \sqrt{\beta }},\frac{2 m_2+m_3}{2 \sqrt{\beta }},-\frac{\sqrt{3} m_3}{2 \sqrt{\beta }},-\frac{2 m_1+m_3}{2 \sqrt{\beta }},\frac{\sqrt{3} \left(m_1+m_2\right)}{2 \sqrt{\beta }},-\frac{m_1-m_2}{2 \sqrt{\beta }})^\top.
\end{equation}

Then
\begin{equation}\label{eigenvectors}
\begin{array}{c}
  \mathcal{E}_3=\mathbf{r}_0,\\
  \mathcal{E}_4 = \mathbf{i}\mathbf{r}_0,\\
  \lambda=\beta ^{3/2},
\end{array}
\end{equation}
and the matrix $\lambda \mathbb{I} + \mathfrak{M}^{-1} \mathfrak{B}$ is
\begin{displaymath}
\beta ^{3/2}\left(
\begin{array}{cccccc}
 \frac{4 m_1+9 m_3}{4}  & -\frac{3 \sqrt{3} m_3}{4}  & m_2 & 0 & -\frac{5 m_3}{4} & \frac{3 \sqrt{3} m_3}{4} \\
 -\frac{3 \sqrt{3} m_3}{4}  & \frac{8 m_2-m_3+4}{4}  & 0 & -2 m_2 & \frac{3 \sqrt{3} m_3}{4} & \frac{m_3}{4} \\
 m_1 & 0 & \frac{4 m_2+9 m_3}{4}  & \frac{3 \sqrt{3} m_3}{4} & -\frac{5 m_3}{4}  & -\frac{3 \sqrt{3} m_3}{4}  \\
 0 & -2 m_1 & \frac{3 \sqrt{3} m_3}{4} & \frac{12-8 m_2-9 m_3}{4}  & -\frac{3 \sqrt{3} m_3}{4}  & \frac{m_3}{4} \\
 -\frac{5 m_1}{4}  & \frac{3 \sqrt{3} m_1}{4} & -\frac{5 m_2}{4}  & -\frac{3 \sqrt{3} m_2}{4}  & \frac{9-5 m_3}{4}  & \frac{3\sqrt{3} \left(m_2-m_1\right)}{4}  \\
 \frac{3 \sqrt{3} m_1}{4} & \frac{m_1}{4} & -\frac{3 \sqrt{3} m_2}{4}  & \frac{m_2}{4} & \frac{3\sqrt{3} \left(m_2-m_1\right)}{4}  & \frac{m_3+3}{4} \\
\end{array}
\right)
\end{displaymath}
By (\ref{Hessian111}), a straightforward computation shows that
\begin{equation}\label{eigenvalues}
  \lambda_1=\frac{3}{2}  \left(1-\sqrt{1-3 \beta }\right) \beta ^{3/2}, ~~~~
  \lambda_2=\frac{3}{2} \left(\sqrt{1-3 \beta }+1\right) \beta ^{3/2}. ~~~~\nonumber
\end{equation}
 As a result, (\ref{threebodylinearstability}) becomes
 \begin{equation}
\beta \leq \frac{1}{27},\nonumber
\end{equation}
or more precisely, Lagrange relative equilibrium  is spectrally   stable if and only if
 \begin{equation}\label{linearstabilityn}
m_1 m_2+ m_3 m_2+ m_1 m_3 \leq \frac{(m_1+m_2+m_3)^2}{27}.
\end{equation}
Moreover, it is easy to see that Lagrange relative equilibrium  is linearly   stable if and only if
 \begin{equation}
m_1 m_2+ m_3 m_2+ m_1 m_3 < \frac{(m_1+m_2+m_3)^2}{27}.\nonumber
\end{equation}
This result has been proved by Gascheau in 1843 \cite{gascheau1843examen} and Routh in 1875 \cite{routh1874laplace} respectively.

Without loss of generality, suppose $m_1\geq m_2\geq m_3$. Then it is easy to see that \eqref{linearstabilityn} yields that
 \begin{equation}
m_1>\frac{1}{18} \left(\sqrt{69}+9\right)>0.961478, m_2+m_3<0.038521.\nonumber
\end{equation}

\begin{figure}
  \center
  \includegraphics[width=8cm]{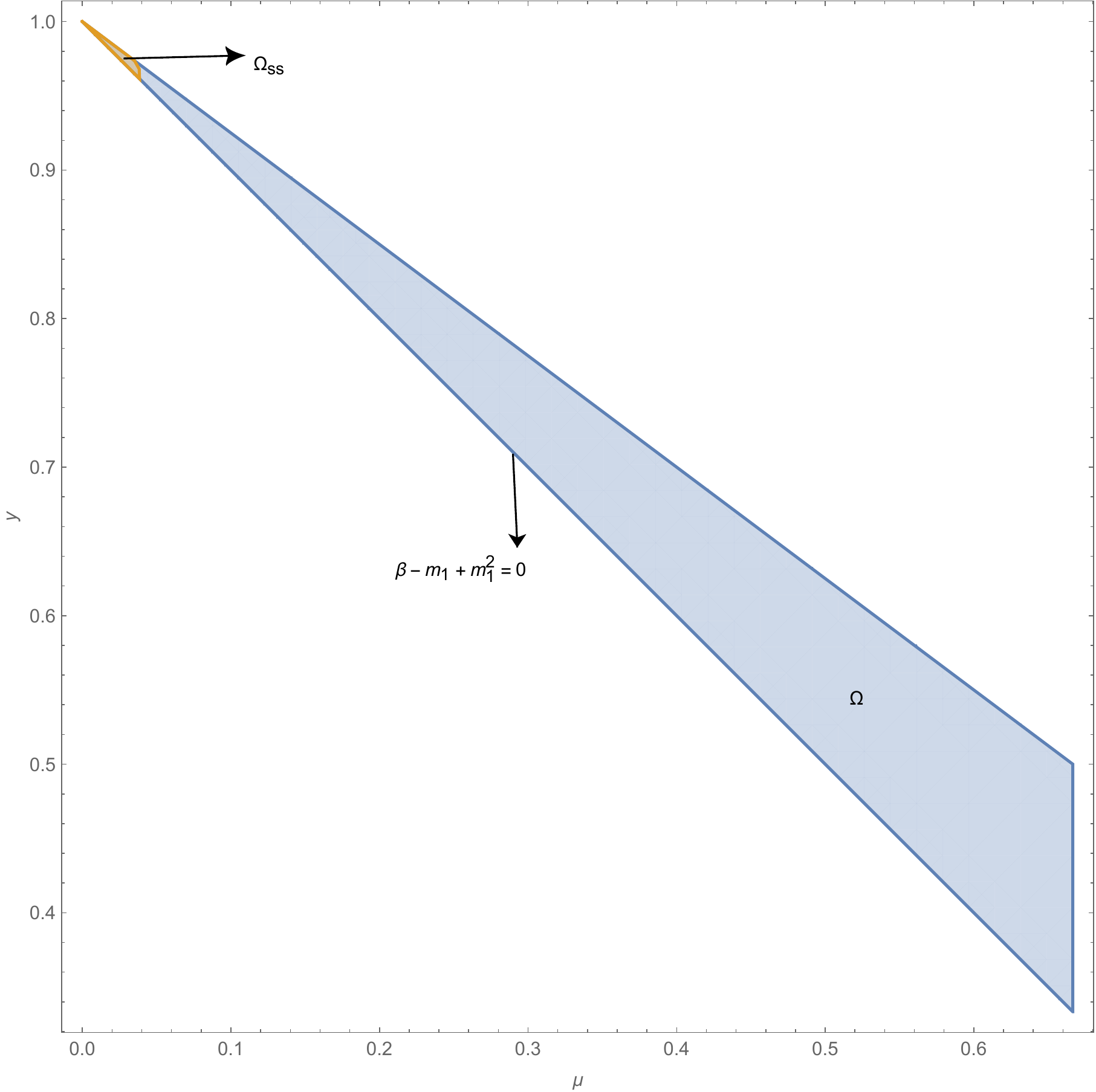}~~~~~\includegraphics[width=8cm]{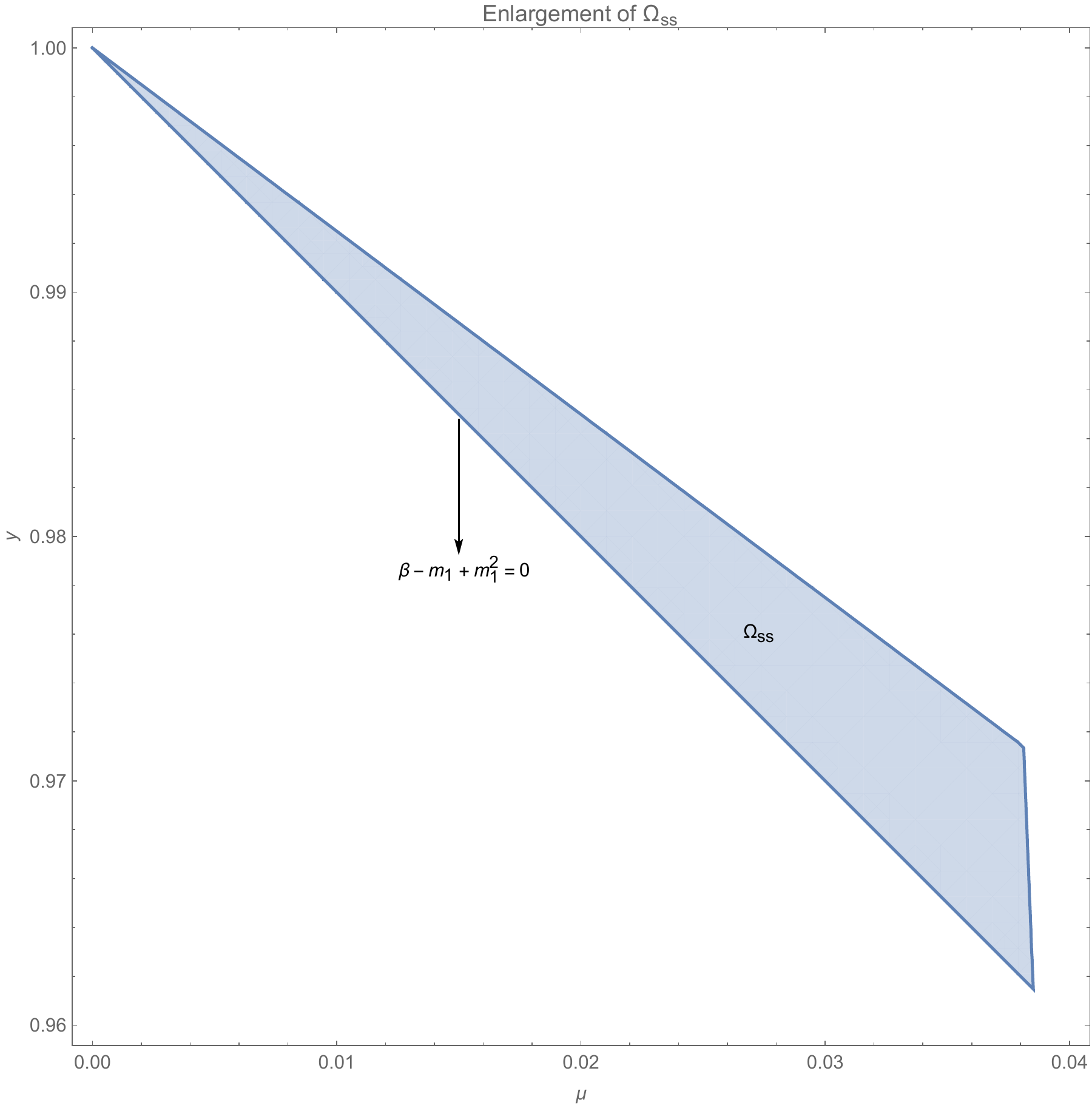}\\
  \caption{spaces $\Omega$ and $\Omega_{ss}$ of masses}\label{spaces of masses}
\end{figure}

Let $\Omega$ be the space of masses of the planar three-body problem, then $\Omega$ could  be represent as
\begin{equation}
\Omega=\{(\beta,m_1): \beta\in(0,\frac{1}{3}], m_1\in[\frac{1}{3},1), \beta- m_1(1-m_1)>0, 4\beta\leq 1+2 m_1-3 m_1^2\}.\nonumber
\end{equation}
In the following it suffices to consider the subset $\Omega_{ss}$ of $\Omega$ corresponding to spectral stability:
\begin{equation}\label{space of masses}
\Omega_{ss}=\{(\beta,m_1)\in \Omega: \beta\in(0,\frac{1}{27}], m_1\in(\frac{\sqrt{69}+9}{18} ,1)\}.
\end{equation}
To make the direct-viewing understanding of sizes of geometric areas $\Omega$ and $\Omega_{ss}$ etc, we would better draw their pictures in a new system of variables $\mu,y$ via the diffeomorphism:
\begin{equation}
\left\{
             \begin{array}{lr}
             \beta = y \mu, & \\
             m_1=1-\mu. &
             \end{array}
\right.\nonumber
\end{equation}
The spaces $\Omega$ and $\Omega_{ss}$  in the variables $\mu,y$ can be seen Figure \ref{spaces of masses}. Obviously,  the space $\Omega_{ss}$ is much smaller than $\Omega$.

Some tedious computation further yields the corresponding eigenvectors
\begin{equation}\label{eigenvectors}
\begin{array}{c}
  \mathcal{E}_1=\left(\frac{m_1-m_3 }{\frac{\kappa m_1}{m_3}},\frac{ 3 m_2-2 \alpha-1}{\frac{\sqrt{3} \kappa m_1}{m_3}},\frac{ m_2-\alpha-m_1}{\frac{ \kappa m_2}{m_3}},\frac{ \alpha+3 m_3-1}{\frac{\sqrt{3} \kappa m_2}{m_3}},\frac{\alpha-m_2+m_3}{\kappa},\frac{\alpha+3 m_1-1}{\sqrt{3}\kappa}\right)^\top, \\[16pt]
  \mathcal{E}_2 = \mathbf{i}\mathcal{E}_1,
\end{array}
\end{equation}
where
\begin{equation}
  \kappa=\sqrt{\frac{4 \beta  m_3 \left(2-6 \beta +\alpha-3 \alpha m_2\right)}{3 m_1 m_2}}, ~~~~\nonumber
\end{equation}
\begin{equation}
\alpha=\sqrt{1-3 \beta }.\nonumber
\end{equation}

\section{Classical Results of Hamiltonian System}
\indent\par

In this section, let us recall some necessary aspects of Hamiltonian system.

We consider an analytic Hamiltonian system, with $n$ degrees of freedom,
having the origin as an equilibrium point:
\begin{equation}\label{Hamiltonian system}
H(p,q) = \sum_{j\geq 2} H_j (p,q),
\end{equation}
where $H_j$ is a homogeneous polynomial of degree $j$ in $(p,q)$ for every $j\geq 2$.

Since  we are interesting in stability, we confine ourself to  the eigenvalues of the
quadratic part $H_2$ of the Hamiltonian are all distinct and purely imaginary.  Then
in suitable symplectic
coordinates, the
quadratic part $H_2$ takes the form
\begin{equation}
H_2 = \sum_{j= 1}^{n} \frac{\omega_j (p_j^2+q_j^2)}{2}.
\end{equation}
Here every $\omega_j$ is
called  a characteristic frequency, and $\varpi=(\omega_1,\cdots,\omega_n)$ is
called  the frequency vector.

\begin{definition}
A frequency vector $\varpi$  satisfies a resonance
relation of order $l>0$ if there exists a linear relationship
\begin{equation}\label{resonant}
 (k, \varpi) = k_1 \omega_1 + \cdots +k_n \omega_n=0,
\end{equation}
where $k=(k_1,   \cdots,  k_n) \in \mathbb{Z}^n$  such that $|k|=|k_1| + \cdots +|k_n| =l$.
\end{definition}

\begin{definition}
A frequency vector $\varpi$  is said to be $(c,\upsilon)$-Diophantine for some $c, \upsilon>0$ if  we have
\begin{equation}
| (k, \varpi)|\geq \frac{c}{|k|^\upsilon}, ~~~~~~~~~\forall k \in \mathbb{Z}^n ~such ~that ~~|k|\neq 0. \nonumber
\end{equation}
A $(c,\upsilon)$-Diophantine frequency vector $\varpi$  is also said to be strongly
incommensurable.
\end{definition}

\begin{definition}
A Birkhoff normal form of degree $l$ for the Hamiltonian (\ref{Hamiltonian system}) is
a polynomial of degree $l$ in symplectic  variables $x,y$ that is actually
a polynomial of degree $[l/2]$ in the variables $\varrho_j=\frac{x_j^2+y_j^2}{2}$.
\end{definition}

Given $l\geq 4$, assume that the frequency vector $\varpi$ is nonresonant up to order $l$.
The well-known Birkhoff theorem \cite{birkhoff1927dynamical} states that, in some
neighbourhood of the origin, there exists a symplectic change of
variables $(p,q)\mapsto (x,y)$,
near to the identity map, such that  in
the new variables the Hamiltonian function is reduced to a Birkhoff normal
form $\mathcal{H}_{Bl}(\varrho)$ of degree $l$ up to terms of degree higher than $l$:
\begin{equation}
H(p,q)=\mathcal{H}(x,y)=\mathcal{H}_{Bl}(\varrho)+ O(|x|+|y|)^{l+1}. \nonumber
\end{equation}

Let us consider a nearly-integrable Hamiltonian written in action-angle
variables $\varrho,\varphi$ defined by $x_j=\sqrt{2\varrho_j}\cos\varphi_j,y_j=\sqrt{2\varrho_j}\sin\varphi_j$:
\begin{equation}\label{Hamiltonian system1}
\mathcal{H}(\varrho,\varphi)=\mathcal{H}_{Bl}(\varrho)+ \mathcal{R}_l(\varrho,\varphi),
\end{equation}
where $\mathcal{R}(\varrho,\varphi)=O(\|\varrho\|)^{[l/2]+1}$, here $\|\varrho\| = \max_{1\leq j\leq n}|\varrho_j|$.

Let us recall the important concepts of non-degenerate and isoenergetically non-degenerate (see \cite{zbMATH05031968}):
\begin{definition}
The Hamiltonian system (\ref{Hamiltonian system}) or (\ref{Hamiltonian system1}) is called to be non-degenerate
in a neighbourhood of the origin if
\begin{equation}
det\left( \frac{\partial^2 \mathcal{H}_{Bl}}{\partial \varrho^2}|_{\varrho=0}\right)\neq 0;\nonumber
\end{equation}
The Hamiltonian system (\ref{Hamiltonian system}) or (\ref{Hamiltonian system1}) is called to be  isoenergetically non-degenerate
in a neighbourhood of the origin if
\begin{equation}
 det\left(
\begin{array}{cc}
 \frac{\partial^2 \mathcal{H}_{Bl}}{\partial \varrho^2}|_{\varrho=0} & \varpi^\top \\
 \varpi  & 0 \\
\end{array}
\right) \neq 0.\nonumber
\end{equation}
\end{definition}

Then it is well known that:
\begin{theorem}\label{KAM} \emph{(\textbf{KAM} \cite{arnol1963small,moser1968lectures})}
In a neighbourhood of an equilibrium point, a non-degenerate or isoenergetically non-degenerate Hamiltonian with a nonresonant frequency vector  up to order $4$  has invariant tori close to the tori of the
linearized system. These tori form a set whose relative measure in the polydisc $\|\varrho\|<\epsilon$ tends to 1 as $\epsilon\rightarrow 0$. In an isoenergetically non-degenerate system
such tori occupy a larger part of each energy level passing near the equilibrium
position.
\end{theorem}
Furthermore, on the relative measure of the set of invariant tori in the polydisc $\|\varrho\|<\epsilon$ we have
\begin{theorem}\label{KAM1}\emph{(\cite{poschel1982integrability,Delshams1996Estimates})}
Consider a non-degenerate or isoenergetically non-degenerate Hamiltonian in a neighbourhood of an equilibrium point. If the frequency vector $\varpi$  is  nonresonant up to order $l\geq 4$, then the relative measure of the set of invariant tori in the polydisc $\|\varrho\|<\epsilon$ is
at least $1-O(\epsilon^{\frac{l-3}{4}})$.  If the frequency vector $\varpi$  satisfies the strong
incommensurability condition, i.e., $(c,\upsilon)$-Diophantine condition, then this measure is $1-O(\exp (-\tilde{c}\epsilon^{\frac{-1}{\upsilon+1}}))$ for a positive number $\tilde{c} = const$.
\end{theorem}

Let us further consider effective stability of a nearly-integrable Hamiltonian
\begin{equation}\label{Hamiltonian system2}
\mathcal{H}(\varrho,\varphi)=(\varrho,\varpi)+ \mathcal{R}_2(\varrho,\varphi).
\end{equation}

First, under $(c,\upsilon)$-Diophantine condition we have
\begin{theorem}\label{effective stability1}\emph{(\cite{giorgilli1998problem,giorgilli1988rigorous,Giorgilli1989Effective})}
Assume the frequency vector $\varpi$  is $(c,\upsilon)$-Diophantine, then  there exist  $c_1,c_2, c_3 = const>0$ such that for every orbit $(\varrho(t),\varphi(t))$ of (\ref{Hamiltonian system2}), with $\|\varrho(0)\|< \epsilon$, one
has
\begin{equation}
\|\varrho(t)-\varrho(0)\|\leq c_1\epsilon^{3}   ~~~~~~~~~~~~~~for~~~ |t|\leq c_3 \exp(c_2\epsilon^{\frac{-1}{\upsilon+1}}),\nonumber
\end{equation}
provided $\epsilon$ is sufficiently small.
\end{theorem}

Although all of the $(c,\upsilon)$-Diophantine frequency vectors
are abundant in measure, however, non-Diophantine
frequency vectors form a dense open set in the space of frequency vectors. Therefore, Diophantine frequency vectors could be
 quite exceptional in some sense.

\begin{definition}\emph{(\cite{niederman2006hamiltonian})}
Let $h$ be a  real analytic in the vicinity of the closed ball $\overline{B_r}$ of radius $r > 0$ in $\mathbb{R}^n$ and has no critical points in $\overline{B_r}$. Then $h$
is steep if and only if its restriction $h|_{\mathcal{P}}$ to any proper affine subspace $\mathcal{P}\subset \mathbb{R}^n$ admits only
isolated critical points.
\end{definition}

\begin{definition}
Let $h$ be a polynomial of degree $l$ in $\varrho_1, \cdots, \varrho_n$ such that
\begin{equation}
h(\varrho)= h_1(\varrho)+\cdots h_l(\varrho),\nonumber
\end{equation}
where is $h_j(\varrho)$ a homogeneous polynomial of degree $k$ in $\varrho_1, \cdots, \varrho_n$. We say that the function $h$ is
\begin{itemize}
  \item convex at $\varrho = 0$, if the quadratic form $h_2(\varrho)$ is either positive or negative definite;
  \item quasi-convex at $\varrho = 0$, if
  \begin{equation}
 h_1(\varrho)=0, h_2(\varrho)=0 ~~~~~~~~~\Rightarrow ~~~~~~~~~\varrho = 0;\nonumber
\end{equation}
  \item directionally quasi-convex at $\varrho = 0$, if
  \begin{equation}
 h_1(\varrho)=0, h_2(\varrho)=0,  \varrho_1, \cdots, \varrho_n \geq 0   ~~~~~~\Rightarrow ~~~~~~\varrho = 0.\nonumber
\end{equation}
\end{itemize}
\end{definition}
\begin{theorem}\label{Nekhoroshev}\emph{(\cite{benettin1998nekhoroshev,fasso1998nekhoroshev,niederman1998nonlinear,poschel1999nekhoroshev})}
Consider the Hamiltonian (\ref{Hamiltonian system1}) in a neighbourhood of the origin $\varrho = 0$. Assume the frequency vector $\varpi$ is nonresonant up to order $l\geq 4$ and the unperturbed Hamiltonian $\mathcal{H}_{Bl}$ is a (directionally) quasi-convex function, then  there exist two positive
constants $a,b$  such that, for sufficiently small $\epsilon$, any orbit $(\varrho(t),\varphi(t))$ of (\ref{Hamiltonian system1}), with $\|\varrho(0)\|< \epsilon$,
satisfies
\begin{equation}
\|\varrho(t)-\varrho(0)\|\leq c_1\epsilon^{a}   ~~~~~~~~~~~~~~for~~~ |t|\leq c_3 \exp(c_2\epsilon^{-b}),\nonumber
\end{equation}
here $c_1,c_2, c_3= const>0$.
\end{theorem}

\begin{remark}
Two positive
constants $a,b$ in Theorem \ref{Nekhoroshev}  can be chosen as $a=\frac{1+\sigma}{n+\sigma}, b=\frac{1}{n+\sigma}$ for any $\sigma \geq 0$, for instance, $a=\frac{1}{n}, b=\frac{1}{n}$ or $a=\frac{1}{2}, b=\frac{1}{2n}$. It should be no surprise that  two
constants $a,b$ in Theorem \ref{Nekhoroshev} are worse than that in Theorem \ref{effective stability1}.
\end{remark}
\begin{remark}
Note that there are some differences between Theorem \ref{Nekhoroshev} and the celebrated Nekhoroshev theorem \cite{nekhorshev1977exponential}. In his
celebrated 1977 article \cite{nekhorshev1977exponential}, Nekhoroshev conjectured that, if the function $\mathcal{H}_{Bl}$ is steep, a weaker condition than convex or quasi-convex properties, the Theorem \ref{Nekhoroshev}  is also correct, however, this conjecture is not complete answered up to now. Note that directionally quasi-convex function may be not steep.
\end{remark}

\section{The Birkhoff Normal Form}
\indent\par
To discuss the stability of relative equilibria, it would be better to employ Hamiltonian form of the $N$-body problem.

\subsection{The Hamiltonian  Near Relative Equilibria}
\indent\par
It is easy to see that the  system (\ref{equations of motion1}) is essentially a Lagrangian system with the Lagrangian function\begin{equation}
\mathcal{L}_{\mathcal{J}}(z,r,\dot{z},\dot{r})= \frac{\dot{r}^2}{2}+ \frac{(1+r)^2}{2} \left[\dot{z}^2_0+|\dot{z}|^2 -\left(\dot{z}^\top Q {z}- \frac{\mathcal{J}}{(1+r)^2}\right)^2\right]+\frac{U(z)}{1+r}.\nonumber
\end{equation}Recall that we have assumed that $\rho=1$.

It follows from the Legendre Transform that
the corresponding Hamiltonian is
\begin{eqnarray}\label{Hamiltonian}
H(r,z,s,w) & =& s\dot{r}+ w^\top \dot{r} -\mathcal{L}_{\mathcal{J}}(z,r,\dot{z},\dot{r})\nonumber\\
& =&\frac{ s^2}{2 }+ \frac{ (1+r)^2}{2} [\dot{z}^2_0+|\dot{z}|^2 + \frac{ \mathcal{J}^2}{(1+r)^4}-(\dot{z}^\top Q {z})^2] - \frac{ U(z)}{1+r},\nonumber
\end{eqnarray}
where
\begin{equation}
\begin{array}{lr}
s=\frac{\partial {\mathcal{L}_{\mathcal{J}}}}{\partial \dot{r}}= \dot{r},\\
  w_k=\frac{\partial {\mathcal{L}_{\mathcal{J}}}}{\partial \dot{z}_k}= (1+r)^2 \left[\frac{ \dot{z}^\top z}{z^2_0} z_k+ \dot{z}_k +\sum_{j=1}^{2N-4}q_{kj} {z}_j(\frac{ \mathcal{J}}{(1+r)^2}-\dot{z}^\top Q z)\right].
\end{array}
\nonumber
\end{equation}
Moreover, a straightforward computation shows that:
\begin{displaymath}
\begin{array}{c}
  \dot{z}^2_0 = \frac{ (\dot{z}^\top z)^2}{z^2_0} =(\frac{ z^\top w}{(1+r)^2} )^2 + O(\parallel (z,w)\parallel^6),\\[8pt]
  \dot{z}_k= \frac{w_k-\mathcal{J}\sum_{j=1}^{2N-4}q_{kj} {z}_j}{(1+r)^2} -\frac{ \dot{z}^\top z}{z^2_0}z_k +\dot{z}^\top Q {z}\sum_{j=1}^{2N-4}q_{kj} {z}_j\\
  =\frac{w_k-\mathcal{J}\sum_{j=1}^{2N-4}q_{kj} {z}_j}{(1+r)^2}- \frac{ z^\top w}{ (1+r)^2}z_k+\sum_{j=1}^{2N-4}q_{kj} {z}_j(\frac{w^\top Q {z}-\mathcal{J}{|z|}^2}{(1+r)^2})+O(\parallel (z,w)\parallel^5).
\end{array}
\end{displaymath}

As a result, we have
\begin{equation}\label{Hamiltonian1}
    \begin{array}{l}
      H(r,z,s,w) =\frac{ s^2}{2}+ \frac{1}{2(1+r)^2} [\mathcal{J}^2+\sum_{k=1}^{2N-4}(w_k-\mathcal{J}\sum_{j=1}^{2N-4}q_{kj} {z}_j)^2 +(w^\top Q {z}-\mathcal{J}{|z|}^2)^2\\[6pt]
      ~~~~~~~~~~~~~~~~~~-~ (z^\top w)^2+ O(\parallel (z,w)\parallel^6)]- \frac{ U(z)}{1+r}\\[10pt]
      =\frac{\omega^2}{2}-\omega^2 r+\frac{ s^2+3\omega^2 r^2+|w|^2-2\omega w^\top Q {z}+ \omega^2 |z|^2}{2}-{r(2 \omega^2 r^2+|w|^2-2\omega w^\top Q {z}+ \omega^2 |z|^2)} \\[6pt]
     +\frac{5\omega^2 r^4+ 3r^2(|w|^2-2\omega w^\top Q {z}+ \omega^2 |z|^2)}{2} + \frac{(w^\top Q {z}-\mathcal{J}{|z|}^2)^2-(z^\top w)^2}{2}+ O(\parallel (r,z,w)\parallel^5)- \frac{ U(z)}{1+r}.
    \end{array}
\end{equation}
We remark that the relative equilibrium $\rho  e^{\mathbf{i}\omega t}\mathbf{r}_0$ is just reduced to the origin, an equilibrium point, of an analytic Hamiltonian system with Hamiltonian (\ref{Hamiltonian1}).

\subsection{The Hamiltonian  Near Lagrange Triangular Point}
\indent\par
For the three-body problem, let us further compute  the Hamiltonian $H(r,z,s,w)$ near Lagrange triangular
point $\mathbf{r}_0$ in \eqref{Lagrangetriangularpoint}.
As a matter of notational convenience, set
\begin{displaymath}
q_0=r,q_1=z_1,q_2=z_2,p_0=s,p_1=w_1,p_2=w_2.
\end{displaymath}

First, by \eqref{eigenvectors}, some tedious computation  yields that
\begin{equation}
\begin{aligned}
&U(z)=\mathcal{U}(z_0 \hat{\mathbf{r}}_0+{z_1 \hat{\mathcal{E}}_1}+{z_2 \hat{\mathcal{E}}_2} )=\widetilde{\mathcal{U}}(z_0 \hat{\mathbf{r}}_0+{z_1 \hat{\mathcal{E}}_1}+{z_2 \hat{\mathcal{E}}_2})\\
&=\lambda+\frac{1}{2} ( \lambda _1 q_1^2+ \lambda _2 q_2^2)+a_{30} q_1^3+a_{12} q_1 q_2^2+a_{21} q_1^2 q_2 + a_{03} q_2^3\\
&+ a_{40} q_1^4 +a_{13} q_1 q_2^3+a_{22} q_1^2 q_2^2+a_{31} q_1^3 q_2+a_{04} q_2^4 +\cdots,
\end{aligned}\nonumber
\end{equation}
where $\cdots$ denotes the terms of degree higher than 4, and
\begin{equation}
\begin{array}{lr}
  a_{30}= (13 \alpha +5)  a_{3012},\\
  a_{12}=-3 (9 \alpha +5)   a_{3012},\\
a_{3012}=\frac{m_3 \beta ^3 (\alpha -1)^2   \left(\alpha -3 m_1+1\right)\left[3 \alpha  (2 \alpha -1)+(2-5 \alpha ) \left(m_2-m_3\right)+3 m_1 \left(3 \alpha -2 m_2+2 m_3\right)\right] }{36 \sqrt{3} \kappa ^3 m_1^2 m_2^2};
\end{array}
\nonumber
\end{equation}
\begin{equation}
\begin{array}{lr}
  a_{21}= -3 (9 \alpha -5) a_{2103}, \\
  a_{03}= (13 \alpha -5) a_{2103}, \\
  a_{2103} = -\frac{m_3 \beta ^3(\alpha +1)^2   \left(\alpha +3 m_1-1\right) \left[10 \alpha ^2+\alpha -9 \alpha  m_2+9 \alpha  m_3-3 (\alpha -4) m_1-18 m_1^2-2\right]}{108 \kappa ^3 m_1^2 m_2^2};
\end{array}
\nonumber
\end{equation}
\begin{equation}
\begin{array}{lr}
  a_{40}= \left[-\left(\alpha ^2-1\right)^2 \left(40 \alpha ^3-123 \alpha ^2+35\right)-27 \left(41 \alpha ^4+8 \alpha ^3+50 \alpha ^2-35\right) m_1 m_2 m_3\right]a_{4004}, \\
  a_{04}= \left[\left(\alpha ^2-1\right)^2 \left(40 \alpha ^3+123 \alpha ^2-35\right)-27 \left(41 \alpha ^4-8 \alpha ^3+50 \alpha ^2-35\right) m_1 m_2 m_3\right] a_{4004},\\
  a_{4004}=\frac{ \left[m_1 \left(8-13 \alpha ^2-4 \alpha +3 (4 \alpha -1) \left(m_2-m_3\right)\right)+(\alpha -1) \left(1-8 \alpha ^2+\alpha +(7 \alpha -1) \left(m_2-m_3\right)\right)+3 (4 \alpha -7) m_1^2+18 m_1^3\right]}{1296 \kappa ^4 m_1^3 m_2^4/(\beta ^{7/2} m_3)};
\end{array}\nonumber
\end{equation}
\begin{equation}
\begin{array}{c}
 a_{22}=\frac{a_{4004}\left[(297 \beta -64) +\frac{3 \left(297 \beta ^2-527 \beta +96\right) m_1}{\beta }+\frac{\left(-1485 \beta ^2+2504 \beta -448\right)   m_1^2}{\beta ^2}+\frac{\left(891 \beta ^2-1284 \beta +224\right) m_1^3 \left(2 \beta +\left(1-m_1\right){}^2\right)}{\beta ^3}\right]}{m_2 m_3/(54\beta ^3)} ;
\end{array}
\nonumber
\end{equation}
\begin{equation}
\begin{array}{lr}
  a_{13}=-a_{31}
  =\frac{ \left[\beta +m_3 (3 m_3-2)\right]\left(\alpha +3 m_1-1\right)\left[8 \alpha ^2-\alpha -7 \alpha  \left(m_2-m_3\right)+m_1 \left(3 \alpha -3 \left(m_2-m_3\right)+6\right)-9 m_1^2+\left(m_2-m_3\right)-1\right]}{4 \sqrt{3} \kappa ^4 m_1^3 m_2^3/(35 \beta ^{11/2} m_3)}.
\end{array}
\nonumber
\end{equation}

It follows that the Hamiltonian \eqref{Hamiltonian1} becomes
\begin{equation}\label{Hamiltonian2}
H = -\frac{\omega_0^2}{2}+ H_2 + H_3+ H_4+\cdots,
\end{equation}
where
\begin{equation}
\begin{array}{lr}
 \omega_0=\omega= \beta ^{3/4}, \\
 H_2 = \frac{1}{2}[p_0^2+ p_1^2+ p_2^2+2\omega _0 \left(p_1 q_2-p_2 q_1\right)+\omega _0^2 (q_0^2 +\frac{ (3 \alpha-1)}{2}  q_1^2-\frac{ (3 \alpha+1)}{2}  q_2^2) ], \\
  \begin{aligned}
&H_3 = 2\omega _0 p_2 q_1 q_0 -2\omega _0 p_1 q_2 q_0 -p_2^2 q_0-p_1^2 q_0-\omega _0^2 q_0^3 -a_{03} q_2^3\\
&-a_{12} q_1 q_2^2-a_{21} q_1^2 q_2-a_{30} q_1^3-\frac{\omega _0^2}{4} q_0  \left((3 \alpha+1) q_1^2+(1-3 \alpha) q_2^2\right) ,
\end{aligned} \\
  \begin{aligned}
&H_4 = \frac{3}{2} p_1^2 q_0^2+\frac{3}{2} p_2^2 q_0^2-\frac{1}{2} p_1^2 q_1^2+\frac{1}{2} p_2^2 q_1^2+\frac{1}{2} p_1^2 q_2^2-\frac{1}{2} p_2^2 q_2^2-3\omega _0 p_2 q_1 q_0^2 +3\omega _0 p_1 q_2 q_0^2 \\
&-\omega _0 p_2 q_1^3 +\omega _0 p_1 q_2^3 -\omega _0 p_2 q_1 q_2^2 +\omega _0 p_1 q_1^2 q_2 +\omega _0^2 \frac{3}{2} q_0^4-2 p_1 p_2 q_2 q_1 +a_{21} q_0 q_2 q_1^2\\
&+a_{30} q_0 q_1^3-a_{31} q_2 q_1^3+ a_{12} q_0 q_2^3+ a_{03} q_0 q_1 q_2^2-a_{13} q_1 q_2^3+(\frac{\omega _0^2}{2}-a_{04}) q_2^4 +(\frac{\omega _0^2}{2}-a_{40})q_1^4 \\
&+ (\omega _0^2-a_{22})q_1^2 q_2^2 +\frac{3\omega _0^2}{4} q_0^2  \left((\alpha+1) q_1^2-(\alpha-1) q_2^2\right).
\end{aligned}
\end{array}\nonumber
\end{equation}
Note that, without loss of generality, we will sometimes omit  the constant term $-\frac{\omega^2}{2}$ in (\ref{Hamiltonian2}) in the following content.

Our task now is to look for a change of variables from $(p,q)$ to
$(x,y)$ such that $H_2$
 takes the form
\begin{equation}
 \frac{\omega_0(x_0^2+y_0^2)}{2}- \frac{\omega_1(x_1^2+y_1^2)}{2}+\frac{\omega_2(x_2^2+y_2^2)}{2}.\nonumber
\end{equation}

Let $\mathbb{J}$ denote the usual symplectic matrix $\left(
                                                       \begin{array}{cc}
                                                           & -\mathbb{I} \\
                                                         \mathbb{I} &   \\
                                                       \end{array}
                                                     \right)
$.
A straight forward computation shows that  the eigenvalues of the matrix $\mathbb{J}\frac{\partial^2 H_2}{\partial^2 (p,q)}$ are
\begin{equation}\label{eigenvalues1}
  \pm \omega_0 \mathbf{i}, ~~~~\pm \omega_1 \mathbf{i}, ~~~~\pm \omega_2 \mathbf{i}, ~~~~ ~~~~\nonumber
\end{equation}
where 
\begin{equation}
\begin{array}{lr}
  \omega _1=\mu_1 \omega_0,~~~~ ~~~~ ~~~~ \mu_1=\sqrt{\frac{{1-\sqrt{1-27 \beta }} }{{2}}}, \\
  \omega _2=\mu_2 \omega_0,~~~~ ~~~~ ~~~~ \mu_2=\sqrt{\frac{{1+\sqrt{1-27 \beta }} }{{2}}}.
  \end{array} ~~~~\nonumber
\end{equation}
Note that  we can restrict our
attention to the variables $p_1,p_2,q_1,q_2$.
For the  eigenvalues $\omega_1 \mathbf{i},\omega_2 \mathbf{i}$, the corresponding
eigenvectors are
\begin{equation}
             \begin{array}{c}
               \left(\frac{3 \alpha-\gamma}{4} \beta ^{3/4},\frac{\mathbf{i} \left(3 \alpha+\gamma\right) \omega _1}{3 \alpha+\gamma-4},\frac{4 \mathbf{i}  \omega _1 \beta ^{-3/4}}{3 \alpha+\gamma-4},1\right)^\top ~~~\text{and}~~~
            \left(\frac{3 \alpha+\gamma}{4}  \beta ^{3/4},\frac{\mathbf{i} \left(\gamma-3 \alpha\right) \omega _2}{-3 \alpha+\gamma+4},-\frac{4 \mathbf{i}  \omega _2 \beta ^{-3/4}}{-3 \alpha+\gamma+4},1\right)^\top,
             \end{array}
             \nonumber
\end{equation}
where
\begin{equation}
\gamma= \sqrt{1-27 \beta }.\nonumber
\end{equation}

It follows that
we can  introduce the following symplectic transformation to reduce the Hamiltonian:
\begin{equation}\label{symplectic transformation}
\left\{
             \begin{array}{lr}
             p_0  =  \sqrt{\omega}x_0 &  \\
             q_0  =  \frac{y_0}{\sqrt{\omega}} &  \\
             p_1  = \frac{\omega _0 (3 \alpha-\gamma  )}{4  \sqrt{\frac{2\gamma  \omega _1}{4-3 \alpha -\gamma }}}x_1+\frac{\omega _0 (3 \alpha +\gamma )}{4  \sqrt{\frac{2\gamma  \omega _2}{4-3 \alpha +\gamma }}}y_2&  \\
              p_2  =  \frac{(\gamma -3 \alpha ) \sqrt{\frac{\gamma  \omega _2}{-3 \alpha +\gamma +4}}}{\sqrt{2} \gamma }x_2-\frac{(3 \alpha +\gamma ) \sqrt{\frac{\gamma  \omega _1}{-3 \alpha -\gamma +4}}}{\sqrt{2} \gamma }y_1&  \\
              q_1  =  -\frac{2 \sqrt{2} \sqrt{\frac{\gamma  \omega _2}{-3 \alpha +\gamma +4}}}{\gamma  \omega _0}x_2-\frac{2 \sqrt{2} \sqrt{\frac{\gamma  \omega _1}{-3 \alpha -\gamma +4}}}{\gamma  \omega _0}y_1&  \\
              q_2  = \frac{1}{\sqrt{2} \sqrt{\frac{\gamma  \omega _1}{-3 \alpha -\gamma +4}}}x_1+\frac{1}{\sqrt{2} \sqrt{\frac{\gamma  \omega _2}{-3 \alpha +\gamma +4}}}y_2&
             \end{array}
\right.
\end{equation}
In fact, by the   transformation \eqref{symplectic transformation}, it follows that the Hamiltonian \eqref{Hamiltonian2} becomes
\begin{equation}\label{Hamiltoniannewnew}
H(x,y) =\frac{\omega_0 (x_0^2+y_0^2)}{2}- \frac{\omega_1(x_1^2+y_1^2)}{2}+\frac{\omega_2(x_2^2+y_2^2)}{2}+ H_3(x,y)+H_4(x,y) +\cdots.\nonumber
\end{equation}

But we'd better  introduce the following  complex symplectic transformation to reduce the Hamiltonian:

\begin{equation}\label{fusymplectictransformation}
\left\{
             \begin{array}{lr}
             x _0=\frac{\zeta_0}{\sqrt{2}}+\frac{\textbf{i} \eta_0}{\sqrt{2}} &  \\
             y _0=\frac{\eta_0}{\sqrt{2}}+\frac{\textbf{i} \zeta_0}{\sqrt{2}} &  \\
             x _1=\frac{\zeta_1}{\sqrt{2}}+\frac{\textbf{i} \eta_1}{\sqrt{2}}&  \\
              y _1=\frac{\eta_1}{\sqrt{2}}+\frac{\textbf{i} \zeta_1}{\sqrt{2}}&  \\
              x _2=\frac{\zeta_2}{\sqrt{2}}+\frac{\textbf{i} \eta_2}{\sqrt{2}}&  \\
              y _2=\frac{\eta_2}{\sqrt{2}}+\frac{\textbf{i} \zeta_2}{\sqrt{2}}&
             \end{array}
\right.
\end{equation}
Then, by the   transformation \eqref{fusymplectictransformation}, it follows that the Hamiltonian \eqref{Hamiltonian2} becomes
\begin{equation}
H(\zeta,\eta)  = \textbf{i}\omega _0 \zeta _0 \eta _0 - \textbf{i}\omega _1 \zeta _1 \eta _1 +\textbf{i} \omega _2 \zeta _2 \eta _2 + H_3(\zeta,\eta)+H_4(\zeta,\eta) +\cdots.
\nonumber
\end{equation}
Note that a  formal series
\begin{equation}
f=\sum_{k,l}f_{k,l}\zeta^k \eta^l, ~~~~ ~~~~k,l\in \mathbb{N}^3
\nonumber
\end{equation}
in the variables $(\zeta,\eta)\in \mathbb{C}^{6}$ represents a real formal series in the  variables $(x,y)$ if and only if
\begin{equation}
f_{k,l}= \textbf{i}^{|k+l|}\overline{f_{l,k}}.
\nonumber
\end{equation}

We will further perform a change of variables $(\zeta,\eta)\mapsto (u,v)$
with a generating function
\begin{equation}
 u _0 \eta _0 +u _1 \eta _1 +u _2 \eta _2 + S_3(u,\eta)+S_4(u,\eta) +\cdots,
\nonumber
\end{equation}
 such that in the new variables $(u,v)$ the Hamiltonian function reduces to a Birkhoff normal form of degree 4 up to terms of degree higher than 4:
\begin{equation}\label{Birkhoff normal form0}
\begin{aligned}
&H(\zeta,\eta) = H(\zeta(u,v),\eta(u,v)) = \mathcal{H}(u,v) \\
&=\textbf{i} \omega _0 u_0 v_0 - \textbf{i}\omega _1 u_1 v _1 +\textbf{i} \omega _2 u _2 v_2 -\frac{1}{2}[\omega_{00}(u_0 v_0)^2+\omega_{11}(u_1 v_1)^2+\omega_{22}(u_2 v_2)^2\\
&+2\omega_{01}(u_0 v_0 u_1 v_1)+2\omega_{02}(u_0 v_0 u_2 v_2)+2\omega_{12}(u_1 v_1 u_2 v_2)] +\cdots,
\end{aligned}
\end{equation}
 where $S_3$ and $S_4$
are forms of degree 3 and 4 in $u,\eta$,
and
\begin{equation}
\zeta=u+\frac{\partial S_3}{\partial \eta}+\frac{\partial S_4}{\partial \eta}+\cdots,~~~~ ~~~~ ~~~~v=\eta+\frac{\partial S_3}{\partial u}+\frac{\partial S_4}{\partial u}+\cdots.\nonumber
\end{equation}

First of all, it is easy to see that $\beta=\frac{1}{27}$ yields that frequency vector $\varpi=(\omega _0,-\omega _1,\omega _2)$  satisfies a resonance
relations of order 2; and
all of resonance
relations of order 3 or 4 satisfied by $\varpi$ are
\begin{equation}\label{allresonancerelations}
\left\{
             \begin{array}{lr}
             \omega _0-2 \omega _1=0 ~~~~~~iff~~~\beta = \frac{1}{36},&  \\
             \omega _2-2 \omega _1=0 ~~~~~~iff~~~\beta = \frac{16}{675},&  \\
             \omega _0-3 \omega _1=0~~~~~~iff~~~\beta = \frac{32}{2187},&  \\
              \omega _0+\omega _1-2 \omega _2=0~~~~~~iff~~~\beta = \frac{64}{1875},&  \\
              \omega _2-3 \omega _1=0~~~~~~iff~~~\beta = \frac{1}{75}.&
             \end{array}
\right.
\end{equation}
For other values of $\beta$, we make use of the relation
\begin{equation}\label{connectrelation}
H(u+\frac{\partial S_3}{\partial \eta}+\frac{\partial S_4}{\partial \eta}+\cdots,\eta) = \mathcal{H}(u,\eta+\frac{\partial S_3}{\partial u}+\frac{\partial S_4}{\partial u}+\cdots)
\end{equation}
to find the  Birkhoff normal form of degree 4.

Equating  the forms of  order 3 in $u,\eta$ of (\ref{connectrelation}) we obtain
\begin{equation}\label{connectrelation3}
\textbf{i} \omega _0 (\frac{\partial S_3}{\partial \eta_0} \eta _0-\frac{\partial S_3}{\partial u_0} u _0) - \textbf{i}\omega _1 (\frac{\partial S_3}{\partial \eta_1} \eta _1-\frac{\partial S_3}{\partial u_1} u _1) +\textbf{i} \omega _2 (\frac{\partial S_3}{\partial \eta_2} \eta _2-\frac{\partial S_3}{\partial u_2} u _2) + H_3(u,\eta)= 0.
\end{equation}
It follows that $S_3$ can be  determined. 

Then by equating  the forms of  order 4 in $u,\eta$ of (\ref{connectrelation}) we obtain
\begin{equation}\label{connectrelation4}
\begin{aligned}
&\textbf{i} \omega _0 (\frac{\partial S_4}{\partial \eta_0} \eta _0-\frac{\partial S_4}{\partial u_0} u _0) - \textbf{i}\omega _1 (\frac{\partial S_4}{\partial \eta_1} \eta _1-\frac{\partial S_4}{\partial u_1} u _1) +\textbf{i} \omega _2 (\frac{\partial S_4}{\partial \eta_2} \eta _2-\frac{\partial S_4}{\partial u_2} u _2) + H_{3\rightarrow 4}+H_4(u,\eta)\\
&+\frac{1}{2}[\omega_{00}(u_0 v_0)^2+\omega_{11}(u_1 v_1)^2+\omega_{22}(u_2 v_2)^2+2\omega_{01}(u_0 v_0 u_1 v_1)+2\omega_{02}(u_0 v_0 u_2 v_2)\\
&+2\omega_{12}(u_1 v_1 u_2 v_2)]=0 ,
\end{aligned}\nonumber
\end{equation}
where $H_{3\rightarrow 4}$ is the forms of  order 4 of $H_3(u+\frac{\partial S_3}{\partial \eta},\eta)$.
It follows that $S_4$ and the  Birkhoff normal form of degree 4 in \eqref{Birkhoff normal form0} can be  determined.

By switching to action-angle variables, the obtained Birkhoff normal form in \eqref{Birkhoff normal form0} becomes
\begin{equation}\label{Birkhoff normal form}
\begin{aligned}
&\mathcal{H}_{B4}
= \omega _0 \varrho_0 - \omega _1 \varrho_1 + \omega _2 \varrho _2  +\frac{1}{2}[\omega_{00}\varrho_0 ^2+\omega_{11}\varrho_1 ^2+\omega_{22}\varrho_2 ^2\\
&+2\omega_{01}\varrho_0 \varrho_1 +2\omega_{02}\varrho_0  \varrho_2 +2\omega_{12}\varrho_1 \varrho_2 ],
\end{aligned}
\end{equation}
where $\varrho_j=\textbf{i}u_j v_j$ ($j=0,1,2$) are action variables, and
\begin{equation}
\begin{array}{c}
  \omega _{00}=-3, \\
  \omega _{01}=-\frac{\sqrt{\gamma +1} \left(21 \gamma ^3-40 \gamma ^2+15 \gamma +4\right)}{12 \sqrt{6} \sqrt{\beta } \gamma  (2 \gamma -1)}, \\
  \omega _{02}=-\frac{\sqrt{\gamma +1} \left(21 \gamma ^2+19 \gamma -4\right)}{4 \sqrt{2} \gamma  (2 \gamma +1)},
\end{array}\nonumber
\end{equation}
\begin{equation}
\begin{array}{c}
  \omega _{11}=\frac{(\gamma -1) \left(1211 \gamma ^4-1336 \gamma ^3+279 \gamma ^2+158 \gamma -76\right)}{72 \gamma ^2 \left(10 \gamma ^2-11 \gamma +3\right)}-\frac{3 \beta ^3 \left(31 \gamma ^2+286 \gamma -236\right)}{8 (\gamma -1) \gamma ^2 (5 \gamma -3) m_1 m_2 m_3}, \\
  \omega _{22}=-\frac{(\gamma +1) \left(1211 \gamma ^4+1336 \gamma ^3+279 \gamma ^2-158 \gamma -76\right)}{72 \gamma ^2 \left(10 \gamma ^2+11 \gamma +3\right)}-\frac{3 \beta ^3 \left(31 \gamma ^2-286 \gamma -236\right)}{8 \gamma ^2 (\gamma +1) (5 \gamma +3) m_1 m_2 m_3},
\end{array}
\nonumber
\end{equation}
\begin{equation}
\begin{array}{c}
  \begin{aligned}
&\omega _{12}=\frac{\sqrt{3 \beta } }{4 (18225 \beta ^2-1107 \beta +16) m_1 m_2 m_3} [(360855 \beta ^2-32265 \beta +624) m_1^3+\\
&(-360855 \beta ^2+32265 \beta -624) m_1^2+3 \beta  (120285 \beta ^2-10755 \beta +208) m_1-4 \beta ^2 (432 \beta +43)].
\end{aligned}
\end{array}
\nonumber
\end{equation}

\indent\par

We remark that
\begin{theorem}\label{Diophantine frequency}
The set $\Gamma_r$ of $\beta \in (0,\frac{1}{27})$ corresponding to resonant frequency vectors $\varpi=(\omega _0,-\omega _1,\omega _2)$  is  countable and dense. The set $\Gamma_d$ of $\beta \in (0,\frac{1}{27})$ corresponding to $(c,\upsilon)$-Diophantine frequency vectors $\varpi=(\omega _0,-\omega _1,\omega _2)$ is a  set of full measure for $\upsilon>6$.
\end{theorem}
{\bf Proof.} 
First, it is easy to see that $\Gamma_r$ is  countable.

Let us recall that
\begin{equation}
\begin{array}{lr}
   \omega_0 = \beta ^{3/4}, ~~~~\omega _1=\mu_1 \omega_0,~~~~\omega _2=\mu_2 \omega_0,\\
   \mu_1=\sqrt{\frac{{1-\sqrt{1-27 \beta }} }{{2}}}, ~~~~ ~~~~ \mu_2=\sqrt{\frac{{1+\sqrt{1-27 \beta }} }{{2}}}.
  \end{array} ~~~~\nonumber
\end{equation}
Set
\begin{displaymath}
\Gamma_0=\{\varpi=(\omega _0,-\omega _1,\omega _2):\beta \in (0,\frac{1}{27})\},
\end{displaymath}
\begin{displaymath}
\Gamma_1=\{(1,-\mu_1,\mu _2):\beta \in (0,\frac{1}{27})\},
\end{displaymath}
\begin{displaymath}
\Gamma_2=\{(\mu _1,\mu_2):\beta \in (0,\frac{1}{27})\}.
\end{displaymath}
Due to
\begin{displaymath}
\mu _1^2+\mu _2^2=1,
\end{displaymath}
geometrically $\Gamma_1, \Gamma_2$ are  circular arcs. We parameterize $\Gamma_2$ as
\begin{displaymath}
\mu _1 = \cos \vartheta,\mu_2= \sin \vartheta, ~~~~ ~~~~\vartheta \in (\frac{\pi}{4},\frac{\pi}{2}).
\end{displaymath}
Then the following mappings are diffeomorphisms:
\begin{displaymath}
\begin{array}{ccccccccc}
  (0,\frac{1}{27}) & \rightarrow &\Gamma_0 & \rightarrow & \Gamma_1 & \rightarrow & \Gamma_2 & \rightarrow & (\frac{\pi}{4},\frac{\pi}{2}):\\
  \beta & \mapsto &\varpi & \mapsto &(1,-\mu_1,\mu_2) & \mapsto & (\mu_1,\mu_2) & \mapsto & \vartheta.
\end{array}
\end{displaymath}

Let us consider  resonance
relations
\begin{equation}
 k_0  - k_1 \mu_1 +k_2 \mu_2=0. \nonumber
\end{equation}
Even if we consider only the case  $k_0=0$,  it is easy to see that the points in $\Gamma_2$ satisfying the relation
\begin{equation}
 - k_1 \mu_1 +k_2 \mu_2=0 \nonumber
\end{equation}
are dense in $\Gamma_2$. Therefore it follows from the diffeomorphisms above that $\Gamma_r$ is also dense in the interval $(0,\frac{1}{27})$.

For fixed $k=(k_0 , k_1,k_2)\in \mathbb{Z}^3, c$ and $\upsilon$ such that $|k|>0, 0<c<1, \upsilon>6$, let us consider the inequality
\begin{equation}
| k_0  - k_1 \cos \vartheta +k_2 \sin \vartheta|< \frac{c}{|k|^\upsilon}. \nonumber
\end{equation}
First,  it is easy to see that
\begin{equation}
 k_1^2 +k_2^2  >0. \nonumber
\end{equation}
Set
\begin{displaymath}
k_1=-\sqrt{k_1^2 +k_2^2}\sin \vartheta_0, ~~~~ ~~~~k_2=\sqrt{k_1^2 +k_2^2}\cos \vartheta_0.
\end{displaymath}
Then we have
\begin{equation}
\frac{-k_0-\frac{c}{|k|^\upsilon}}{\sqrt{k_1^2 +k_2^2}}< \sin (\vartheta+\vartheta_0)< \frac{-k_0+\frac{c}{|k|^\upsilon}}{\sqrt{k_1^2 +k_2^2}}. \nonumber
\end{equation}
But it is easy to show that, for any $\delta>0$, the measure of the set
\begin{equation}
\{x \in [-2\pi, 2\pi]:y-\delta< \sin x < y+\delta\} \nonumber
\end{equation}
does not exceed $\sigma \sqrt{\delta}$, here $\sigma$ is a constant number. Hence the measure of the set of $\vartheta$ such that
\begin{equation}
| k_0  - k_1 \cos \vartheta +k_2 \sin \vartheta|< \frac{c}{|k|^\upsilon} \nonumber
\end{equation}
does not exceed $\frac{\sigma\sqrt{c}}{\sqrt[4]{k_1^2 +k_2^2}|k|^\frac{\upsilon}{2}}(\leq \frac{\sigma\sqrt{c}}{|k|^\frac{\upsilon}{2}})$.

Since the number of values of $k$ with $|k|=n$ does not exceed $25 n^2$, the measure of the set $\Gamma_{c,\upsilon} \subset (\frac{\pi}{4},\frac{\pi}{2})$ of $\vartheta$ such that
\begin{equation}
| k_0  - k_1 \cos \vartheta +k_2 \sin \vartheta|< \frac{c}{|k|^\upsilon} \nonumber
\end{equation}
for any $|k|>0$ does not exceed
\begin{equation}
\sum^{\infty}_{n=1}\frac{25\sigma\sqrt{c} n^2}{n^\frac{\upsilon}{2}} \leq 25\sigma \sqrt{c} \sigma_{\upsilon},\nonumber
\end{equation}
here the constant $\sigma_{\upsilon}$ depends only on $\upsilon$. As $c\rightarrow 0$, the measure of the set $\Gamma_{c,\upsilon}$ tends to zero.
Therefore it follows from the diffeomorphisms above that the measure of the set $(0,\frac{1}{27})\setminus \Gamma_d$  is zero.

The proof of the theorem is now complete.

$~~~~~~~~~~~~~~~~~~~~~~~~~~~~~~~~~~~~~~~~~~~~~~~~~~~~~~~~~~~~~~~~~~~~~~~~~~~~~~~~~~~~~~~~~~~~~~~~~~~~~~~~~~~~~~~~~~~~~~~~~~~~~~~~~~~~~~~~~~~~~~~~~~~~~~~~~~~~~~~~~~~\Box$\\

We conclude this section with the notation of the following space of masses:\begin{figure}
  \center
  \includegraphics[width=10cm]{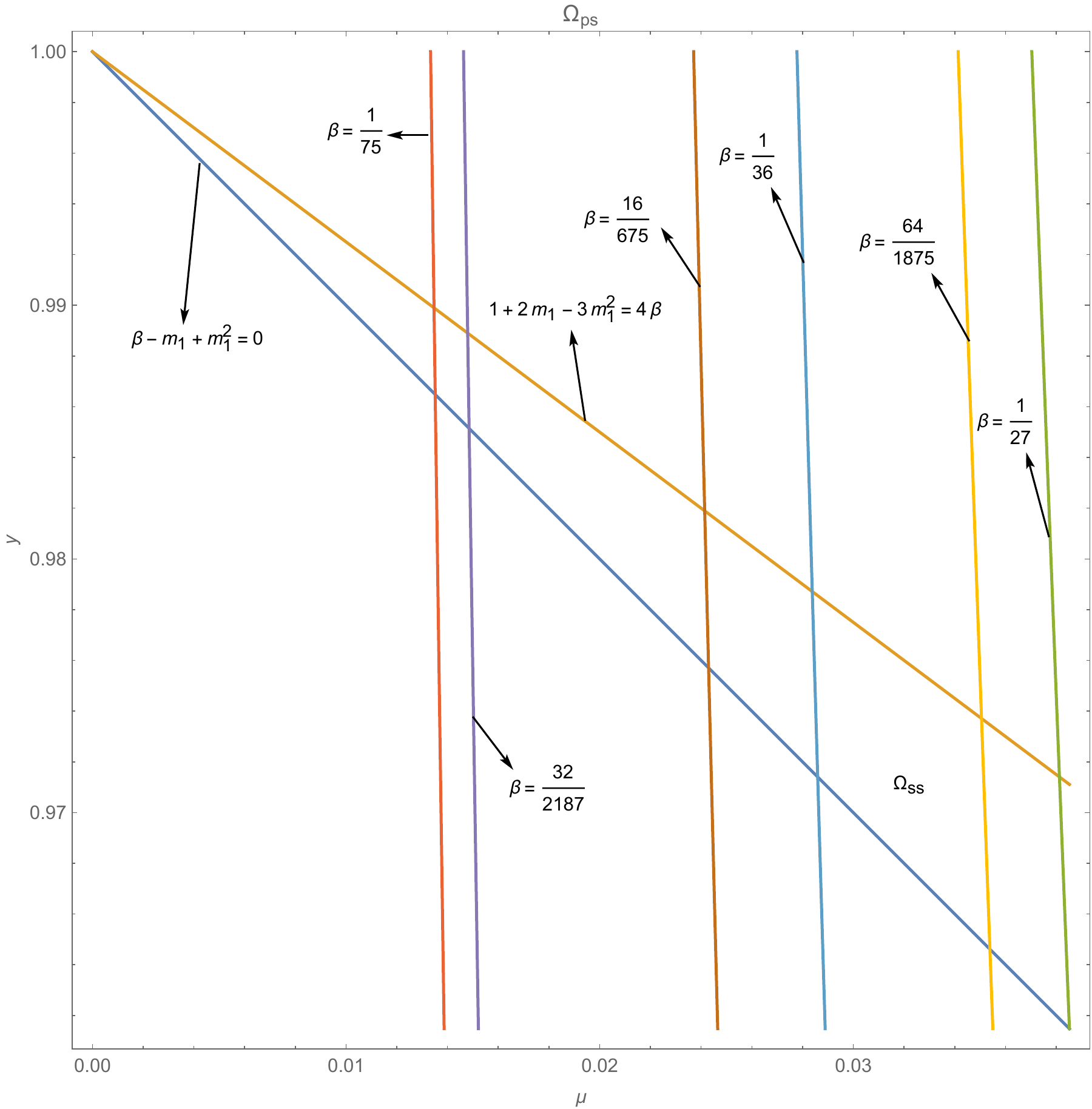}\\
  \caption{masses space $\Omega_{ps}$ and resonance }\label{spaces of massesps}
\end{figure}
\begin{equation}\label{space of massesps}
\Omega_{ps}=\{(\beta,m_1)\in \Omega_{ss}: \beta\in(0,\frac{1}{27})\setminus \{\frac{1}{75},\frac{32}{2187},\frac{16}{675},\frac{1}{36},\frac{64}{1875}\}\},
\end{equation} please see Figure \ref{spaces of massesps}.
Taking the resonance
relations  (\ref{allresonancerelations}) into consideration, from now on, we confine ourselves to the masses in $\Omega_{ps}$.

\section{KAM Stability}
\indent\par

In this section, let us investigate the KAM stability (i.e., stability in the sense of measure) of Lagrange relative equilibrium.  The main result  is the following theorem.
\begin{theorem}\label{KAMapply}
Possibly except the following cases corresponding to resonance
\begin{equation}
   \begin{array}{c}
    \beta=\frac{1}{75}; \beta=\frac{32}{2187}; \beta=\frac{16}{675}; \beta=\frac{1}{36}; \beta=\frac{64}{1875},
   \end{array}\nonumber
\end{equation}
for every choice of   masses  of the planar three-body problem satisfying $\beta<\frac{1}{27}$, there are a great quantity of KAM invariant tori (or quasi-periodic solutions) in a small neighbourhood of Lagrange relative equilibrium. Furthermore, these tori form a set whose relative measure rapidly tends to 1 as the neighbourhood shrinks to zero; in particular, the relative measure exponentially tends to 1 for almost every choice of  the masses.
\end{theorem}
As a corollary, we have the following result:
\begin{theorem}\label{KAMapply2}
Possibly except the following cases corresponding to resonance
\begin{equation}
   \begin{array}{c}
    \beta=\frac{1}{75}; \beta= \frac{32}{2187}; \beta=\frac{16}{675}; \beta=\frac{1}{36}; \beta= \frac{64}{1875}
   \end{array}\nonumber
\end{equation}
for every choice of   masses  of the planar three-body problem satisfying $\beta<\frac{1}{27}$, Lagrange relative equilibrium  is KAM stable.
\end{theorem}

Before proving Theorem \ref{KAMapply}, let us recall a classical result on
real algebraic varieties. 
\begin{definition}\emph{(\cite{whitney1957elementary})}
An algebraic partial manifold $P$ in $\mathbb{R}^n$ is a point set, associated with a number $\nu$,
with the following property. Take any $p\in P$. Then there exists a set of polynomials $f_1, \cdots, f_\nu$ of rank $\nu$ at $p$ (i.e.,  the number of independent
differential $df_1(p), \cdots, df_\nu(p)$ is $\nu$), and a neighborhood $\mathcal{N}$ of $p$, such that $P\bigcap \mathcal{N}$
is the set of zeros in $\mathcal{N}$ of these $f_j$. The number $n-\nu$ is the dimension of the
partial manifold.
\end{definition}
\begin{theorem}\label{real algebraic variety}\emph{(\cite{whitney1957elementary})}
Let $V \subset \mathbb{R}^n$ be a real algebraic variety, then $V$ can be  split as  a union of  a finite number of  partial algebraic
manifolds:
\begin{equation}
V=P_1\bigcup P_2\bigcup \cdots \bigcup P_s, \nonumber
\end{equation}
each $P_j$ being an algebraic partial manifold in $V$, and the $P_j$ being disjoint. Here, the dimension $n_j$ of $P_j$ are decrease. Furthermore, $s\leq 2^n-1$ and each $P_j$ has but a finite number of
topological components.
\end{theorem}

{\bf Proof of Theorem \ref{KAMapply}:}

First, let us investigate degeneracy and isoenergetical degeneracy of the the Hamiltonian (\ref{Birkhoff normal form}).

A straight forward computation shows that
\begin{equation}\label{degerate}
det\left(
\begin{array}{ccc}
 \omega _{00} & \omega _{01} & \omega _{02} \\
 \omega _{01} & \omega _{11} & \omega _{12} \\
 \omega _{02} & \omega _{12} & \omega _{22} \\
\end{array}
\right)=\frac{-27 \beta  }{128 (16-675 \beta )^2 (1-36 \beta )^2 \gamma ^4 m_1^2 m_2^2 m_3^2} f_{deg}
\end{equation}
and
\begin{equation}\label{isodegerate}
det\left(
\begin{array}{cccc}
 \omega _{00} & \omega _{01} & \omega _{02} & \omega _0 \\
 \omega _{01} & \omega _{11} & \omega _{12} & -\omega _1 \\
 \omega _{02} & \omega _{12} & \omega _{22} & \omega _2 \\
 \omega _0 & -\omega _1 & \omega _2 & 0 \\
\end{array}
\right)=\frac{-27 \beta f_{isodeg} }{64 (16-675 \beta )^2 (1-36 \beta )^2 (1-27 \beta )^2 m_1^2 m_2^2 m_3^2} ,
\end{equation}
where
\begin{equation}\label{degeratef}
\begin{aligned}
& f_{deg} = \frac{2(1-36 \beta )^2}{3}  (52542675 \beta ^3+178185258 \beta ^2-9896841 \beta -47632) \beta ^4\\
&-11 (397050199920 \beta ^5-40790893923 \beta ^4+4055047758 \beta ^3-243771759 \beta ^2+6417616 \beta \\
&-59392) \beta ^3 m_1+(5465578392450 \beta ^6+19309935720393 \beta ^5-3995019640449 \beta ^4\\
&+327340481715 \beta ^3-13039336341 \beta ^2+250520816 \beta -1857536) \beta ^2 m_1^2+(2408448\\
&-15298708984020 \beta ^6-29436067209393 \beta ^5+7048034089254 \beta ^4-562788423405 \beta ^3\\
&+20645100208 \beta ^2-359200768 \beta ) \beta  m_1^3+3 [(1821859464150 \beta ^6+4980794507091 \beta ^5\\
&-1182106602432 \beta ^4+94244985459 \beta ^3-3452615664 \beta ^2+59975680 \beta -401408) \\
&m_1^4 (2 \beta +m_1^2-2 m_1+1)],\nonumber
\end{aligned}
\end{equation}
\begin{equation}\label{isodegeratef}
\begin{aligned}
& f_{isodeg}=(1-36 \beta )^2 (52542675 \beta ^3+178185258 \beta ^2-9896841 \beta -47632) \beta ^4\\
&-6 (1114633724580 \beta ^5-129174146793 \beta ^4+12399204438 \beta ^3-701681085 \beta ^2+17908688 \beta\\
&-163328) \beta ^3 m_1+3 (2856548519100 \beta ^6+9467506918989 \beta ^5-1979796586608 \beta ^4\\
&+162904807989 \beta ^3-6502400730 \beta ^2+125103520 \beta -928768) \beta ^2 m_1^2-24 (992795560920 \beta ^6\\
&+1777266330759 \beta ^5-427262272146 \beta ^4+34351179507 \beta ^3-1270282468 \beta ^2+22280704 \beta\\
&-150528) \beta  m_1^3+9 (952182839700 \beta ^6+2412746489943 \beta ^5-573816097674 \beta ^4\\
&+46035466371 \beta ^3-1699679520 \beta ^2+29762048 \beta -200704) m_1^4 (2 \beta +m_1^2-2 m_1+1).\nonumber
\end{aligned}
\end{equation}

\begin{figure}
  \center
  \includegraphics[width=8cm]{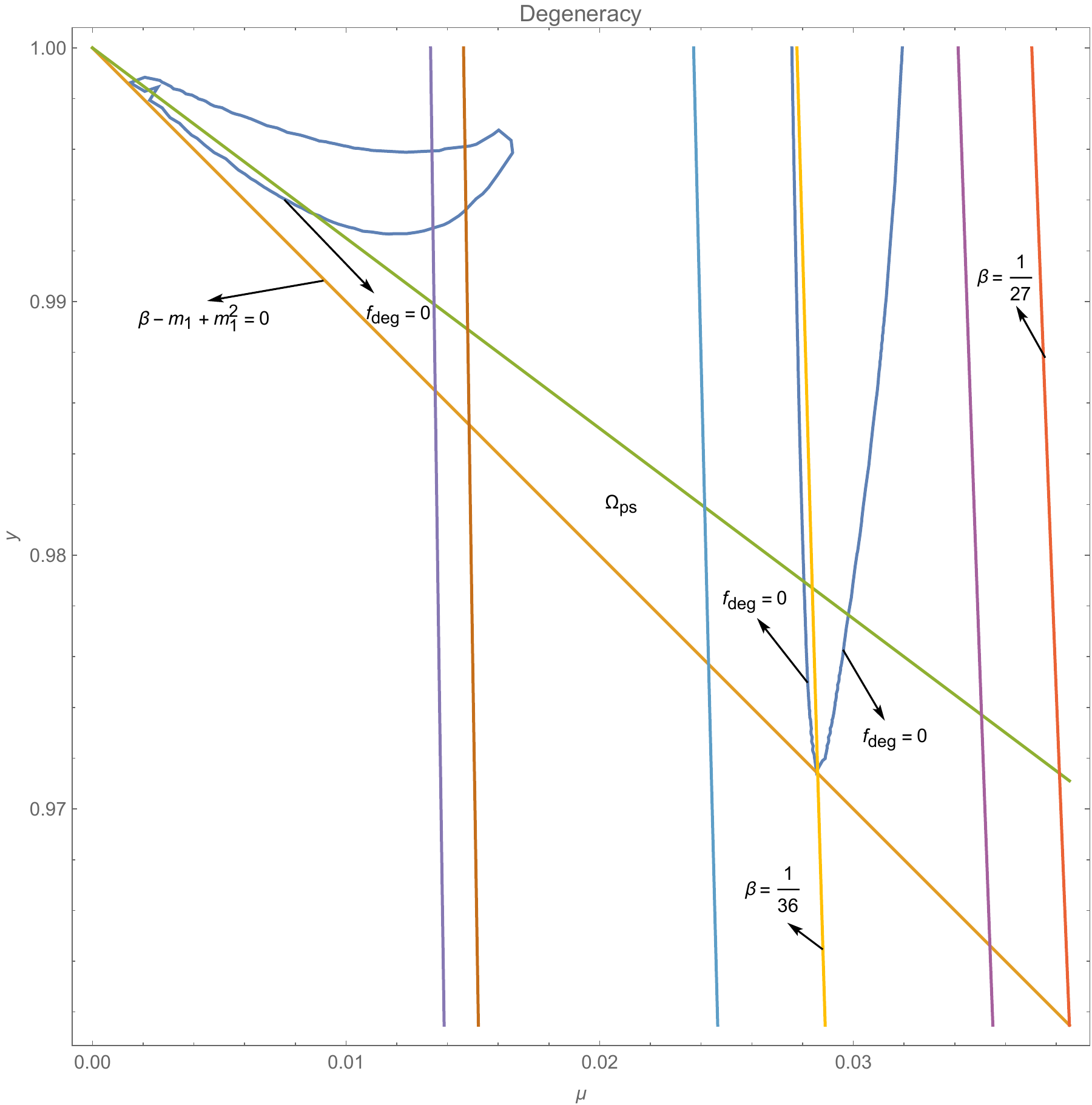}~~~~~\includegraphics[width=8cm]{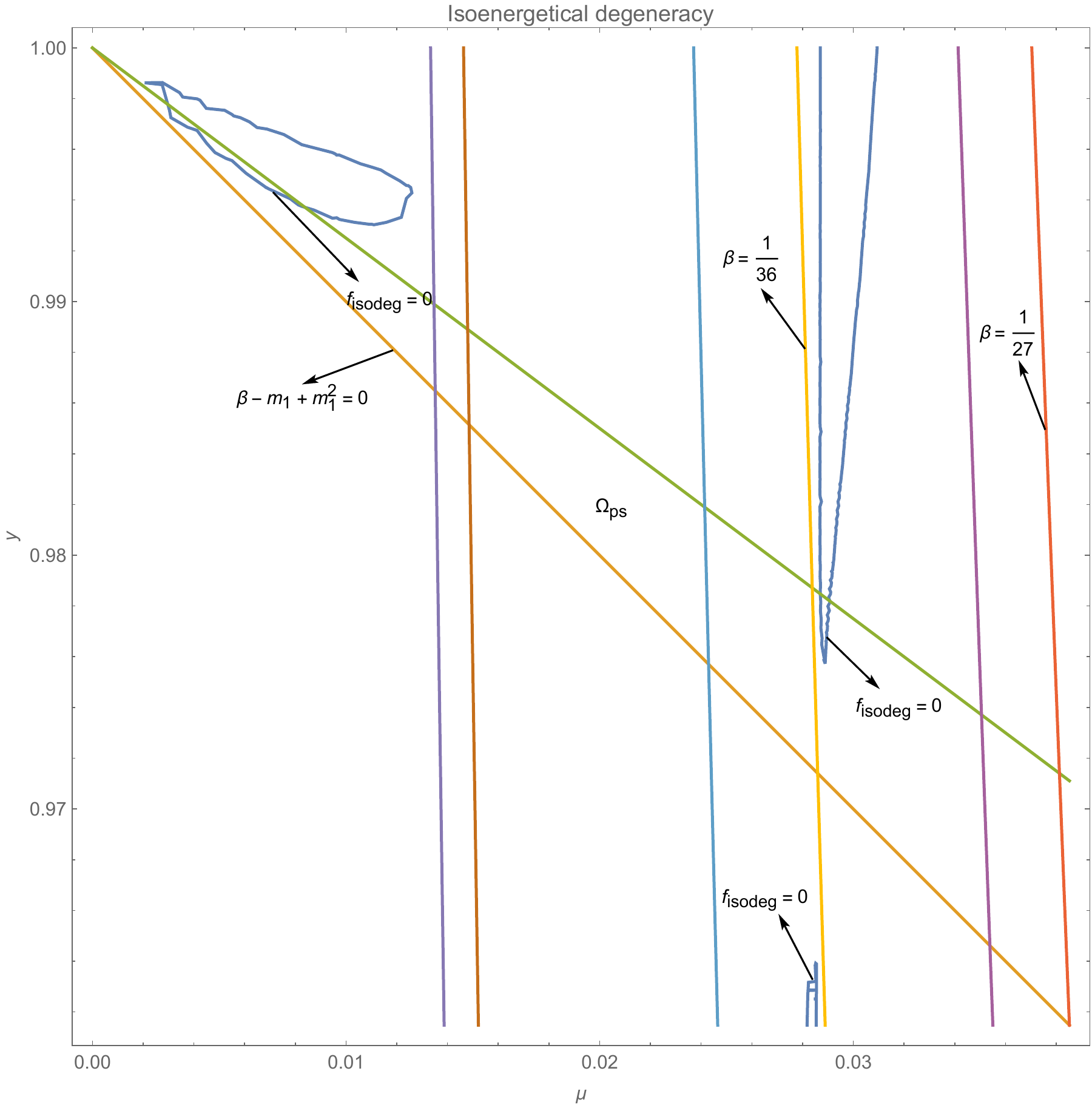}\\
  \caption{plots of $V_f$ and $V_g$} \label{tuihua}
\end{figure}

Therefore, the Hamiltonian (\ref{Birkhoff normal form}) is non-degenerate if and only if $f_{deg}\neq 0$,
and the Hamiltonian (\ref{Birkhoff normal form}) is isoenergetically non-degenerate if and only if $f_{isodeg}\neq 0$.
Thus the set $V_{f_{deg}}$ of points $(\beta,m_1)$ such that the Hamiltonian (\ref{Birkhoff normal form}) is degenerate is  a real algebraic variety. So is the set $V_{f_{isodeg}}$ of isoenergetically degenerate. By Theorem \ref{real algebraic variety},
it follows that each of  $V_{f_{deg}}$ and  $V_{f_{isodeg}}$ is an union of a finite number of zero-dimensional points and one-dimensional ``curves". To make the direct-viewing understanding of the real algebraic varieties  $V_{f_{deg}}$ and  $V_{f_{isodeg}}$, we give the  plots of zero locus sets of $f_{deg}$ and $f_{isodeg}$, please see  Figure \ref{tuihua}.

So the Hamiltonian (\ref{Birkhoff normal form}) is non-degenerate and isoenergetically non-degenerate for almost every choices of $(\beta,m_1)$. Furthermore, a straight forward computation shows that $f$ and $g$ can not be 0  at the same time in the space of masses $\Omega_{ps}$.

As a result, it follows from  Theorem \ref{KAM}, Theorem \ref{KAM1} and  Theorem \ref{Diophantine frequency} that Theorem \ref{KAMapply} holds.

$~~~~~~~~~~~~~~~~~~~~~~~~~~~~~~~~~~~~~~~~~~~~~~~~~~~~~~~~~~~~~~~~~~~~~~~~~~~~~~~~~~~~~~~~~~~~~~~~~~~~~~~~~~~~~~~~~~~~~~~~~~~~~~~~~~~~~~~~~~~~~~~~~~~~~~~~~~~~~~~~~~~\Box$\\

\section{Effective Stability}
\indent\par

We now turn to effective stability.

First, it follows from  Theorem \ref{effective stability1}  and  Theorem \ref{Diophantine frequency} that
\begin{theorem}\label{effective stability1apply}
For almost every choice of  masses of the planar three-body problem satisfying  $\beta\in (0,\frac{1}{27})$, there exists a small neighbourhood $\mathcal{N}$ of Lagrange relative equilibrium  such that, for every orbit $(\varrho(t),\varphi(t))$ whose initial value is in $\mathcal{N}$, one
has
\begin{equation}
\|\varrho(t)-\varrho(0)\|\leq c_1\epsilon^{3}   ~~~~~~~~~~~~~~for~~~ |t|\leq c_3 \exp(c_2\epsilon^{-b}),\nonumber
\end{equation}
provided $\epsilon$ is sufficiently small, here $c_1,c_2, c_3= const>0$ and the constant $b$ is any number in the interval $(0, \frac{1}{7})$. Therefore, Lagrange relative equilibrium  is exponentially  stable for almost every choice of  masses of the planar three-body problem such that $\beta<\frac{1}{27}$.
\end{theorem}

On the one hand, the masses in  Theorem \ref{effective stability1apply} yielding exponential stability 
are abundant in measure. On the other hand,  the masses in  Theorem \ref{effective stability1apply} may be quite exceptional in some sense. For example, according to Theorem \ref{Diophantine frequency},  the masses which can not yield exponential stability are dense. Therefore, it is not allowed to 
have measuring error of masses for applying Theorem \ref{effective stability1apply}, however it is impossible for no measuring error of masses   in practice.

So let us investigate directional quasi-convexity of the Birkhoff normal form $\mathcal{H}_{B4}(\varrho)$. 

As a matter of convenience, first of all, let us give the following result.
\begin{lemma}\label{lemma}
The Birkhoff normal form $\mathcal{H}_{B4}(\varrho)$  is
\begin{itemize}
  \item convex at $\varrho = 0$, if and only if
  \begin{equation}
 det\left(
\begin{array}{cc}
 \omega _{00} & \omega _{01}  \\
 \omega _{01} & \omega _{11}  \\
\end{array}
\right)>0, det\left(
\begin{array}{ccc}
 \omega _{00} & \omega _{01} & \omega _{02} \\
 \omega _{01} & \omega _{11} & \omega _{12} \\
 \omega _{02} & \omega _{12} & \omega _{22} \\
\end{array}
\right) < 0;\nonumber
\end{equation}
  \item quasi-convex at $\varrho = 0$, if and only if
  \begin{equation}
 a_0 a_2-a_1^2 >0 ;\nonumber
\end{equation}

  \item directionally quasi-convex at $\varrho = 0$, if and only if
\begin{equation}
\left\{
             \begin{array}{lr}
             a_0 \neq 0& \\
             a_2=0&\\
             a_0(a_0+2a_1 \frac{\mu_1}{\mu_2})>0 &
             \end{array}
\right.\nonumber
\end{equation}
or
\begin{equation}
\left\{
             \begin{array}{lr}
             a_0 \neq 0 ~and~ a_2\neq 0& \\
             a_0 h(\frac{\mu_1}{\mu_2})>0&\\
             - \frac{a_1}{a_2}) \notin [0, \frac{\mu_1}{\mu_2}]  &
             \end{array}
\right.\nonumber
\end{equation}
or
\begin{equation}
\left\{
             \begin{array}{lr}
             a_0 \neq 0 ~and~ a_2\neq 0& \\
             a_0 h(\frac{\mu_1}{\mu_2})>0&\\
             a_0 h(- \frac{a_1}{a_2}))>0&\\
             - \frac{a_1}{a_2} \in [0, \frac{\mu_1}{\mu_2}]  &
             \end{array}
\right.\nonumber
\end{equation}
\end{itemize}
where
\begin{equation}
 \begin{array}{c}
   a_0= \omega _{00} \mu_1^2+ 2\omega _{01} \mu_1+\omega _{11},\\
   a_1= -\omega _{00} \mu_1\mu_2 + \omega _{02} \mu_1-\omega _{01} \mu_2+\omega _{12}, \\
   a_2= \omega _{00} \mu_2^2- 2\omega _{02} \mu_2+\omega _{22}, \\
   h(x)=a_0+ 2a_1 x + a_2 x^2.
 \end{array}
   \nonumber
\end{equation}
\end{lemma}
{\bf Proof.} 
Recall that
\begin{equation}
\begin{aligned}
&\mathcal{H}_{B4}(\varrho)
= \omega _0 \varrho_0 - \omega _1 \varrho_1 + \omega _2 \varrho _2  +\frac{1}{2}[\omega_{00}\varrho_0 ^2+\omega_{11}\varrho_1 ^2+\omega_{22}\varrho_2 ^2\\
&+2\omega_{01}\varrho_0 \varrho_1 +2\omega_{02}\varrho_0  \varrho_2 +2\omega_{12}\varrho_1 \varrho_2 ].
\end{aligned}
\end{equation}

By $\omega_{00}=-3$, it follows that $\mathcal{H}_{B4}(\varrho)$  is convex at $\varrho = 0$ if and only if the quadratic form
\begin{displaymath}
\omega_{00}\varrho_0 ^2+\omega_{11}\varrho_1 ^2+\omega_{22}\varrho_2 ^2+2\omega_{01}\varrho_0 \varrho_1 +2\omega_{02}\varrho_0  \varrho_2 +2\omega_{12}\varrho_1 \varrho_2
\end{displaymath}
is negative definite, that is,
\begin{equation}
 det\left(
\begin{array}{cc}
 \omega _{00} & \omega _{01}  \\
 \omega _{01} & \omega _{11}  \\
\end{array}
\right)>0, det\left(
\begin{array}{ccc}
 \omega _{00} & \omega _{01} & \omega _{02} \\
 \omega _{01} & \omega _{11} & \omega _{12} \\
 \omega _{02} & \omega _{12} & \omega _{22} \\
\end{array}
\right) < 0.\nonumber
\end{equation}

Thanks to
\begin{displaymath}
\omega _0 \varrho_0 - \omega _1 \varrho_1 + \omega _2 \varrho _2=0,
\end{displaymath}
or
\begin{displaymath}
 \varrho_0 = \lambda _1 \varrho_1 - \lambda _2 \varrho _2,
\end{displaymath}
the quadratic form
\begin{displaymath}
\omega_{00}\varrho_0 ^2+\omega_{11}\varrho_1 ^2+\omega_{22}\varrho_2 ^2+2\omega_{01}\varrho_0 \varrho_1 +2\omega_{02}\varrho_0  \varrho_2 +2\omega_{12}\varrho_1 \varrho_2
\end{displaymath}
reduces to
the quadratic form
\begin{displaymath}
a_0\varrho_1 ^2+a_2\varrho_2 ^2 +2a_1\varrho_1 \varrho_2.
\end{displaymath}
Then it is evident to see that $\mathcal{H}_{B4}(\varrho)$  is quasi-convex at $\varrho = 0$, if and only if
  \begin{equation}
 a_0 a_2-a_1^2 >0;\nonumber
\end{equation}
 $\mathcal{H}_{B4}(\varrho)$  is directionally quasi-convex at $\varrho = 0$, if and only if
\begin{equation}
\left\{
             \begin{array}{lr}
             \lambda _1 \varrho_1 \geq \lambda _2 \varrho _2\geq 0&\\
             a_0\varrho_1 ^2+a_2\varrho_2 ^2 +2a_1\varrho_1 \varrho_2=0 &
             \end{array}
\right.   \Rightarrow \varrho_1 =\varrho _2=0,\nonumber
\end{equation}
or $a_0\neq 0$ and the equation
\begin{displaymath}
h(x)=a_0+ 2a_1 x + a_2 x^2=0
\end{displaymath}
has no roots in the interval  $[0, \frac{\mu_1}{\mu_2}]$.
As a result, it is easy to see that  the theorem holds.

$~~~~~~~~~~~~~~~~~~~~~~~~~~~~~~~~~~~~~~~~~~~~~~~~~~~~~~~~~~~~~~~~~~~~~~~~~~~~~~~~~~~~~~~~~~~~~~~~~~~~~~~~~~~~~~~~~~~~~~~~~~~~~~~~~~~~~~~~~~~~~~~~~~~~~~~~~~~~~~~~~~~\Box$

Let $\Omega_{c},\Omega_{qc},\Omega_{dqc}$ be the subsets of  the space $\Omega_{ps}$ of masses corresponding to convexity, quasi-convexity and directional quasi-convexity respectively. Then a straightforward computation shows that $\Omega_{c}$ is empty, that is,  $\mathcal{H}_{B4}(\varrho)$  is not convex at $\varrho = 0$ for any choice of masses of the three-body problem.

\begin{figure}
  \center
  \includegraphics[width=8cm]{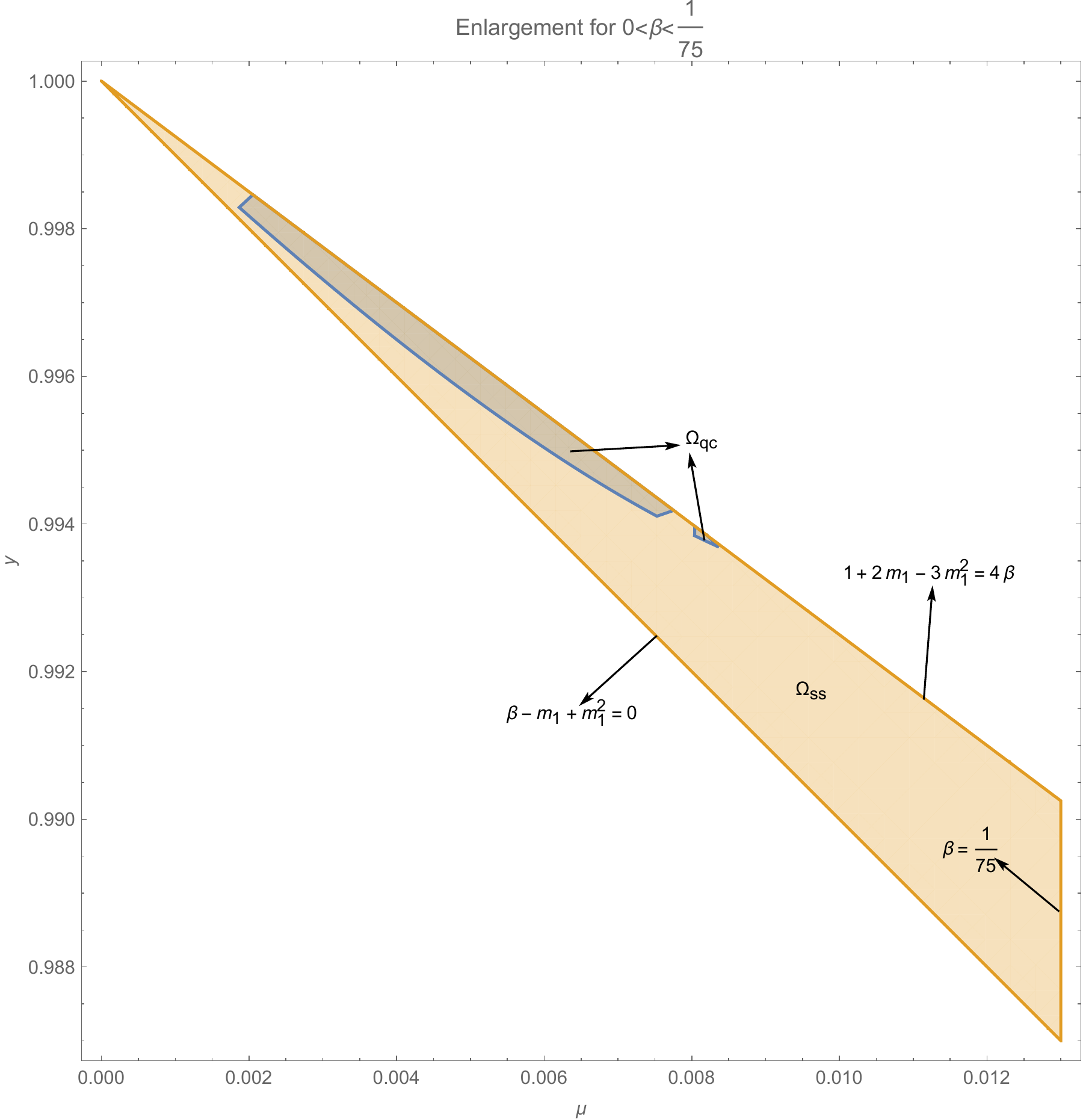}~~~~~\includegraphics[width=8cm]{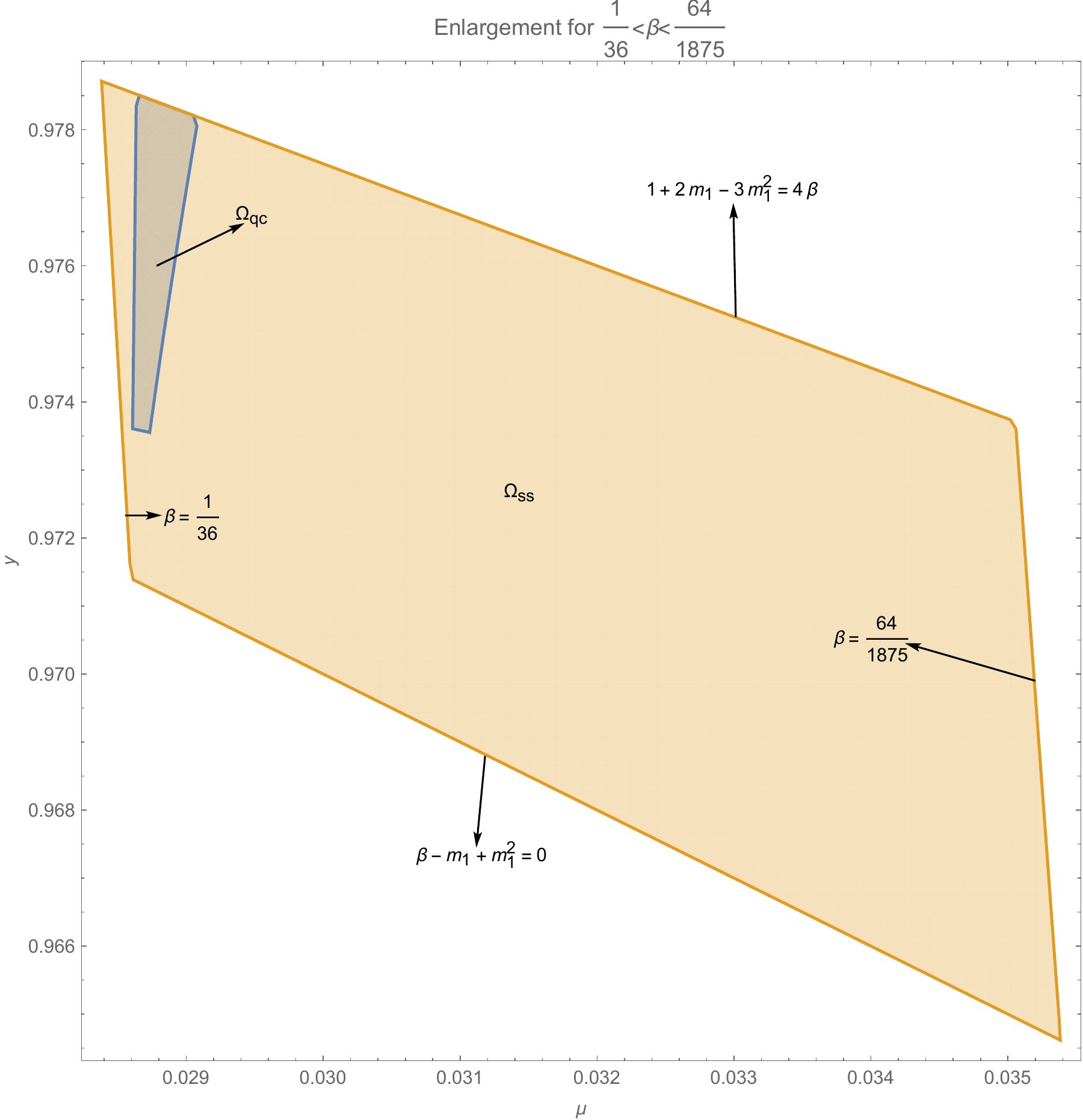}\\
  \caption{enlargement of $\Omega_{qc}$} \label{quasi-convex12}
\end{figure}

For quasi-convexity, a straightforward computation shows that
\begin{equation}
a_0 a_2-a_1^2 = - det\left(
\begin{array}{cccc}
 \omega _{00} & \omega _{01} & \omega _{02} & \omega _0 \\
 \omega _{01} & \omega _{11} & \omega _{12} & -\omega _1 \\
 \omega _{02} & \omega _{12} & \omega _{22} & \omega _2 \\
 \omega _0 & -\omega _1 & \omega _2 & 0 \\
\end{array}
\right),\nonumber
\end{equation}
as a result,
\begin{equation}
a_0 a_2-a_1^2 >0  \Leftrightarrow f_{isodeg}>0,\nonumber
\end{equation}
and $\Omega_{qc}$ is empty for $\frac{1}{75}< \beta <\frac{1}{36}$ and $\frac{64}{1875}< \beta <\frac{1}{27}$ but not for $0< \beta <\frac{1}{75}$ or $\frac{1}{36}< \beta <\frac{64}{1875}$.
To make the direct-viewing understanding of the space $\Omega_{qc}$,  please see  Figure \ref{tuihua} and \ref{quasi-convex12}, note that the picture for $\Omega_{qc}$ is enlarged.

\begin{figure}
  \center
  \includegraphics[width=8cm]{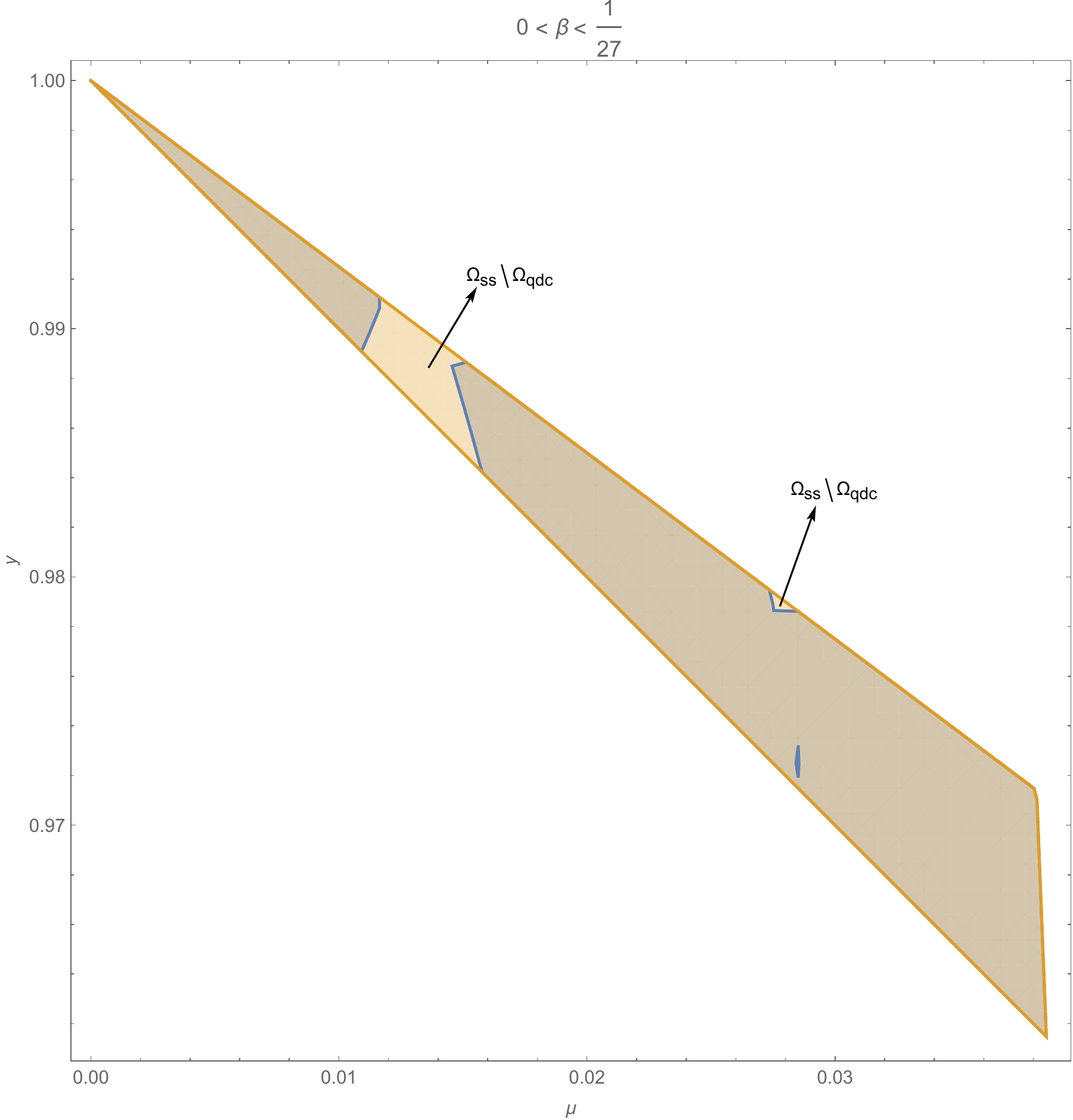}~~~~~\includegraphics[width=8cm]{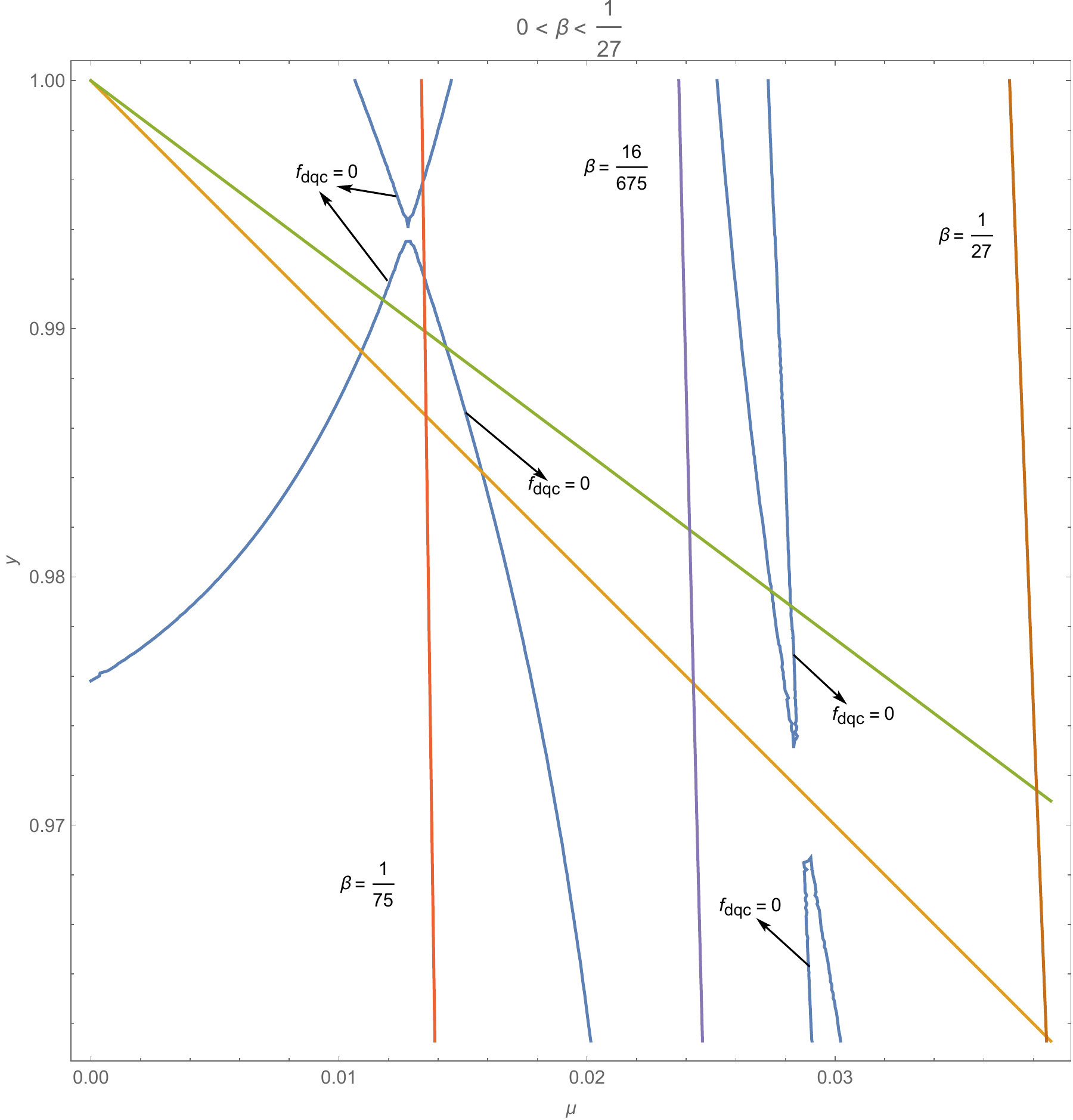}\\
  \caption{plots of $\Omega_{dqc}$} \label{directionally-quasi-convex}
\end{figure}

For directional quasi-convexity, it follows from Lemma \ref{lemma} that
\begin{equation}
\Omega_{dqc}=\Omega_{pdqc} \setminus \Omega_{ndqc},\nonumber
\end{equation}
here
\begin{equation}
\begin{array}{c}
  \Omega_{pdqc}=\{(\beta,m_1)\in \Omega_{ps}: a_0 h(\frac{\mu_1}{\mu_2})>0\}, \\
  \Omega_{ndqc}=\{(\beta,m_1)\in \Omega_{pdqc}:  a_2\neq 0, 0<- \frac{a_1}{a_2}<\frac{\mu_1}{\mu_2}, a_0 h(- \frac{a_1}{a_2})\leq 0\}.
\end{array}
\nonumber
\end{equation}
Some tedious computation  shows that $\Omega_{ndqc}$ is empty, and for $(\beta,m_1)\in \Omega_{ps}$,
\begin{equation}
a_0 h(\frac{\mu_1}{\mu_2})>0 \Leftrightarrow f_{dqc}>0,
\nonumber
\end{equation}
where
\begin{equation}
\begin{array}{lc}
 \begin{aligned}
& f_{dqc} = [\frac{\left(236-62 \gamma ^4-479 \gamma ^3+1299 \gamma ^2-994 \gamma \right) (\gamma +1)^2}{m_1 m_2 m_3}+729 \left(76-401 \gamma ^3+81 \gamma ^2-18 \gamma \right) ]\\
&[\frac{3 \beta  \left(16-469476 \beta ^3+71469 \beta ^2-2199 \beta \right)}{m_1 m_2 m_3}+\left(1509030 \beta ^3+2316519 \beta ^2-133983 \beta +1936\right) ].
\end{aligned}\nonumber
\end{array}
\end{equation}
As a result, the space $\Omega_{dqc}=\Omega_{pdqc}$ is a subset of $\Omega_{ps}$ satisfying $f_{dqc}>0$. It is easy to see that the space $\Omega_{dqc}$ is a large part of the space $\Omega_{ps}$ geometrically. To make the direct-viewing understanding of the space $\Omega_{dqc}$,  please see  Figure \ref{directionally-quasi-convex}.

To sum up, we have
\begin{theorem}\label{directionally quasi-convex}
The Birkhoff normal form $\mathcal{H}_{B4}(\varrho)$  is
\begin{itemize}
  \item never  convex at $\varrho = 0$, i.e., $\Omega_{c}$ is empty;
  \item quasi-convex at $\varrho = 0$, if and only if $(\beta,m_1)\in \Omega_{qc}$, here
  \begin{equation}
 \Omega_{qc}=\{(\beta,m_1)\in \Omega_{ps}: f_{isodeg}>0,~ \beta\in (0,\frac{1}{75})\bigcup (\frac{1}{36},\frac{64}{1875}) \};\nonumber
\end{equation}

  \item directionally quasi-convex at $\varrho = 0$, if and only if $(\beta,m_1)\in \Omega_{dqc}$, here
\begin{equation}
 \Omega_{dqc}=\{(\beta,m_1)\in \Omega_{ps}: f_{dqc} >0\}.\nonumber
\end{equation}
\end{itemize}
\end{theorem}

To study exponential stability, it suffices to consider directional quasi-convexity. Indeed, it follows from Theorem \ref{Nekhoroshev} that
\begin{theorem}\label{Nekhoroshevapply}
For every choice of masses of the planar three-body problem satisfying  $(\beta,m_1)\in \Omega_{dqc}$, there exists a small neighbourhood $\mathcal{N}$ of Lagrange relative equilibrium  such that, for every orbit $(\varrho(t),\varphi(t))$ whose initial value is  in $\mathcal{N}$, one
has
\begin{equation}
\|\varrho(t)-\varrho(0)\|\leq c_1\epsilon^{a}   ~~~~~~~~~~~~~~for~~~ |t|\leq c_3 \exp(c_2\epsilon^{-b}),\nonumber
\end{equation}
provided $\epsilon$ is sufficiently small, here $c_1,c_2, c_3= const>0$ and the constants $a, b$ can be chosen as $a=\frac{1+\sigma}{n+\sigma}, b=\frac{1}{n+\sigma}$ for any $\sigma \geq 0$. Therefore, Lagrange relative equilibrium  is exponentially  stable for  every choice of  masses in the space $\Omega_{dqc}$.
\end{theorem}
\begin{remark}
For a choice of  masses in the space $\Omega_{ps}\setminus \Omega_{dqc}$, one should extend and apply a more general Theorem \ref{Nekhoroshev} to explore exponential  stability. As a matter of fact, there is a more general Theorem \ref{Nekhoroshev} in \cite{benettin1998nekhoroshev} with  a weaker condition than  directional quasi-convexity, the condition is nearer to steepness, that is,  the Birkhoff normal form $\mathcal{H}_{B6}(\varrho)$  is $3$-jet nondegenerate at $\varrho = 0$. But the computation is too complicated to obtain the Birkhoff normal form $\mathcal{H}_{B6}(\varrho)$.
\end{remark}

Let us give the following result that the Birkhoff normal form $\mathcal{H}_{B4}(\varrho)$ is steep, although there is no proof of a general Theorem \ref{Nekhoroshev} under the condition of  steepness at present.

\begin{theorem}\label{steep}
Possibly except the following cases corresponding to resonance
\begin{equation}
   \begin{array}{c}
    \beta=\frac{1}{75}; \beta=\frac{32}{2187}; \beta=\frac{16}{675}; \beta=\frac{1}{36}; \beta=\frac{64}{1875},
   \end{array}\nonumber
\end{equation}
for every choice of   masses  of the planar three-body problem satisfying $\beta<\frac{1}{27}$, $\mathcal{H}_{B4}(\varrho)$ is steep in some neighbourhood ${B_{r_{\beta}}}$ of the origin, where
\begin{displaymath}
r_{\beta}=\frac{\beta^{3/4}}{\sqrt{5+\frac{ 9}{ (1-36 \beta )^2 ( 1-27 \beta)}+\frac{10^6}{ (27 \beta -1)^2 (36 \beta -1)^2 (675 \beta -16)^2 m_1^2 m_2^2 m_3^2}} }.
\end{displaymath}
\end{theorem}
{\bf Proof.} 
We divide our proof in two steps.

First, for every choice of positive masses of the planar three-body problem satisfying  $(\beta,m_1)\in \Omega_{ps}$,  we prove that $\mathcal{H}_{B4}(\varrho)$ has no critical points in $\overline{B_{r_{\beta}}}$. (step 1)

Assume  the point $\varrho=(\varrho_0,\varrho_1,\varrho_2)^\top$ is a critical point of $\mathcal{H}_{B4}(\varrho)$, then
\begin{equation}
\frac{\partial\mathcal{H}_{B4}}{\partial \varrho}=0,\nonumber
\end{equation}
or
\begin{equation}\label{critical point}
W \varrho=-\varpi^\top,
\end{equation}
here
\begin{equation}
W=\frac{\partial^2\mathcal{H}_{B4}}{\partial \varrho^2}=\left(
\begin{array}{ccc}
 \omega _{00} & \omega _{01} & \omega _{02} \\
 \omega _{01} & \omega _{11} & \omega _{12} \\
 \omega _{02} & \omega _{12} & \omega _{22} \\
\end{array}
\right),\nonumber
\end{equation}
and $\varpi=(\omega _0,-\omega _1,\omega _2)$ is the frequency vector.

It follows  that
\begin{equation}\label{critical point1}
\begin{aligned}
\omega_0^2 + \omega _1 ^2 + \omega _2^2 \leq  (\omega_{00}^2+\omega_{11}^2+\omega_{22}^2+2\omega_{01}^2 +2\omega_{02}^2 +2\omega_{12}^2)(\varrho_0 ^2+\varrho_1 ^2+\varrho_2 ^2).
\end{aligned}
\end{equation}

The aim is to prove that (\ref{critical point1}) yields the following estimation
\begin{equation}\label{estimation}
\begin{array}{c}
  \frac{\omega_0^2 + \omega _1 ^2 + \omega _2^2}{\omega_{00}^2+\omega_{11}^2+\omega_{22}^2+2\omega_{01}^2 +2\omega_{02}^2 +2\omega_{12}^2}>  r_{\beta}.
\end{array}
\end{equation}

A straightforward computation shows that the three equations
\begin{displaymath}
\omega_{11}=0, ~~~~~~~~\omega_{22}=0,~~~~~~~~\omega_{12}=0
\end{displaymath}
have no solution for $0<\beta <\frac{1}{27}$ and $0.96<m_1<1$, thus
\begin{displaymath}
\omega_{11}^2+\omega_{22}^2 +2\omega_{12}^2>0.
\end{displaymath}

By virtue of
\begin{equation}
\begin{aligned}
&\omega_{11}^2+\omega_{22}^2 +2\omega_{12}^2=\frac{9}{32 (27 \beta -1)^2 (36 \beta -1)^2 (675 \beta -16)^2 m_1^2 m_2^2 m_3^2} \\
&[\frac{\beta ^4 (1-36 \beta )^2}{3} (52542675 \beta ^4-254067678 \beta ^3+24800847 \beta ^2-387536 \beta +3456)\\
&+6(1092450983580 \beta ^5-177644083959 \beta ^4+12794359374 \beta ^3-478172835 \beta ^2\\
&+8933708 \beta -65728) \beta ^4 m_1  -\beta ^2 m_1^2(8796996017100 \beta ^7+1736483408451 \beta ^6\\
&-260078284944 \beta ^5+7108908219 \beta ^4+541625274 \beta ^3-40867544 \beta ^2+981888 \beta\\
&-8192 ) +4 (6037174483920 \beta ^7-2675577372453 \beta ^6+422084648466 \beta ^5\\
&-35545883265 \beta ^4+1718731704 \beta ^3-47333488 \beta ^2+688128 \beta -4096) \beta  m_1^3\\
&-(8796996017100 \beta ^7-4818222493029 \beta ^6+805786218810 \beta ^5-69657248025 \beta ^4\\
&+3410662284 \beta ^3-94469792 \beta ^2+1376256 \beta -8192) m_1^4 \left(2 \beta +m_1^2-2 m_1+1\right)],
\end{aligned}\nonumber
\end{equation}
it is easy to see that
\begin{displaymath}
\omega_{11}^2+\omega_{22}^2 +2\omega_{12}^2<\frac{10^6}{ (27 \beta -1)^2 (36 \beta -1)^2 (675 \beta -16)^2 m_1^2 m_2^2 m_3^2}
\end{displaymath}
 for $0<\beta <\frac{1}{27}$ and $0.96<m_1<1$.

Thanks to
\begin{displaymath}
2\omega_{01}^2 +2\omega_{02}^2 = \frac{9 \left(775656 \beta ^3-63135 \beta ^2+1701 \beta -16\right)}{8 (1-36 \beta )^2 (27 \beta -1)}
\end{displaymath}
and
\begin{displaymath}
-16<775656 \beta ^3-63135 \beta ^2+1701 \beta -16<-\frac{16}{81} ~~~~~~~~for ~~~~0<\beta <\frac{1}{27},
\end{displaymath}
we have the following estimation
\begin{equation}
\begin{aligned}
&9+\frac{ 2}{9 (1-36 \beta )^2 ( 1-27 \beta)}<\omega_{00}^2+\omega_{11}^2+\omega_{22}^2+2\omega_{01}^2 +2\omega_{02}^2 +2\omega_{12}^2\\
&<9+\frac{ 18}{ (1-36 \beta )^2 ( 1-27 \beta)}+\frac{10^6}{ (27 \beta -1)^2 (36 \beta -1)^2 (675 \beta -16)^2 m_1^2 m_2^2 m_3^2}.
\end{aligned}\nonumber
\end{equation}

Consequently,
\begin{equation}
\begin{array}{c}
  \frac{\omega_0^2 + \omega _1 ^2 + \omega _2^2}{\omega_{00}^2+\omega_{11}^2+\omega_{22}^2+2\omega_{01}^2 +2\omega_{02}^2 +2\omega_{12}^2}>   \frac{\beta^{3/2}}{5+\frac{ 9}{ (1-36 \beta )^2 ( 1-27 \beta)}+\frac{10^6}{ (27 \beta -1)^2 (36 \beta -1)^2 (675 \beta -16)^2 m_1^2 m_2^2 m_3^2} },
\end{array}\nonumber
\end{equation}
and it follows that $\mathcal{H}_{B4}(\varrho)$ has no critical points in $\overline{B_{r_{\beta}}}$.

Our task now is to prove that  a restriction $\mathcal{H}_{B4}|_{\mathcal{P}}$ of $\mathcal{H}_{B4}(\varrho)$ to any proper affine subspace $\mathcal{P}\subset \mathbb{R}^3$ admits only
isolated critical points. (step 2)

Suppose that $\mathcal{H}_{B4}(\varrho)$ admits nonisolated critical points in some proper affine subspace $\mathcal{P}\subset \mathbb{R}^3$,  then by virtue of  (\ref{critical point}) and  $\varpi\neq 0$, it follows  that $\mathcal{P}$ is a two-dimensional plane. Let $\mathcal{P}= Span(\xi_1,\xi_2)$, $\xi_1,\xi_2$ are independent. Then we have the relations
\begin{equation}\label{wuqiongduo}
Rank(W\xi_1, W\xi_2,-\varpi^\top)=Rank(W\xi_1, W\xi_2)=1
\end{equation}
So  the matrix $W$ is noninvertible, and the frequency vector $\varpi^\top$ belongs to the image set of linear operator $W$, in other words,
\begin{equation}\label{wuqiongduo1}
\varpi^\top \in Span(\left(
                            \begin{array}{c}
                              \omega _{00} \\
                              \omega _{01} \\
                              \omega _{02} \\
                            \end{array}
                          \right)
,\left(
                            \begin{array}{c}
                              \omega _{01} \\
                              \omega _{11} \\
                              \omega _{12} \\
                            \end{array}
                          \right),\left(
                            \begin{array}{c}
                              \omega _{02} \\
                              \omega _{12} \\
                              \omega _{22} \\
                            \end{array}
                          \right) ).
\end{equation}

A straight forward computation shows that  for every choice of   masses  satisfying  $(\beta,m_1)\in \Omega_{ps}$, $(\omega_{02},\omega_{12},\omega_{22})$ and $(\omega_{01},\omega_{11},\omega_{12})$ are independent. Hence
\begin{equation}
det\left(
\begin{array}{ccc}
 \omega _{00} & \omega _{01} & \omega _{02} \\
 \omega _{01} & \omega _{11} & \omega _{12} \\
 \omega _{02} & \omega _{12} & \omega _{22} \\
\end{array}
\right)=0, ~~~~~ det\left(
\begin{array}{ccc}
 1 & \omega _{01} & \omega _{02} \\
 -\lambda _{1} & \omega _{11} &  \omega _{12}\\
 \lambda _{2} & \omega _{12} &  \omega _{22}\\
\end{array}
\right)=0.\nonumber
\end{equation}
Some tedious computation  shows that for $(\beta,m_1)\in \Omega_{ps}$, the equations above cannot hold at the same time.

So $\mathcal{H}_{B4}(\varrho)$ is steep in some neighbourhood ${B_{r_{\beta}}}$ of the original point for any $(\beta,m_1)\in \Omega_{ps}$.
The proof of Theorem \ref{steep} is now complete.

$~~~~~~~~~~~~~~~~~~~~~~~~~~~~~~~~~~~~~~~~~~~~~~~~~~~~~~~~~~~~~~~~~~~~~~~~~~~~~~~~~~~~~~~~~~~~~~~~~~~~~~~~~~~~~~~~~~~~~~~~~~~~~~~~~~~~~~~~~~~~~~~~~~~~~~~~~~~~~~~~~~~\Box$

\section{Conclusion}
\indent\par

For the planar three-body problem, based on the moving coordinates introduced in  \cite{yu2019problem}, which allows us to obtain a reduced system of   equations of motion suitable for describing  the motion of  particles near relative equilibria,   we mainly discussed the nonlinear stability  of  Lagrange relative equilibrium.

First,  we proved that  a relative equilibrium   is orbitally stable    if and only if the origin  of the reduced system is Lyapunov stable. Before discussing the nonlinear stability, we gave some well known information on linear stability, and it is clear that it is more convenient to get the information by using the method based on the moving coordinates.

Next, it is necessary to get the  Birkhoff normal
form of the Hamiltonian near Lagrange triangular
point. Although the construction of the normal form is  simple in concept,  but it is difficult to  obtain the normal form. Thus this paper requires some computer assistance. Certainly, intensive computation cannot be avoided in celestial mechanics.

By virtue of the  celebrated KAM theorem,  we proved  that Lagrange relative equilibrium  is KAM stable, except possibly six special resonant cases, if it is spectrally stable.  Indeed,  there are a great quantity of KAM invariant tori in a small neighbourhood of Lagrange relative equilibrium, provided that the mass parameter $\beta\in(0,\frac{1}{27}]$, except possibly six special resonant cases $\beta=\frac{1}{75},\frac{32}{2187}, \frac{16}{675},\frac{1}{36},\frac{64}{1875},\frac{1}{27}$.  Furthermore, these tori or quasi-periodic solutions  form a set whose relative measure rapidly tends to 1.

We also investigated  the effective (exponential) stability of Lagrange relative equilibrium  by  the  celebrated Nekhoroshev's
theory.  First, we proved that Lagrange relative equilibrium  is exponentially stable for almost every choice of positive masses of the planar three-body problem, except a dense but zero measure  set of masses,   if it is spectrally stable. Then we proved that Lagrange relative equilibrium  is exponential stable for any choice of positive masses in a large open subset  of spectrally stable space of masses. This large open subset is described by  directional quasi-convexity.

Finally, we proved that the  Birkhoff normal
form of the Hamiltonian near Lagrange triangular
point is steep provided that the mass parameter $\beta\in(0,\frac{1}{27}]$, except possibly six special resonant cases. This may be useful for further research of Lagrange relative equilibrium.

We hope to further explore  the nonlinear stability problem  of  general relative equilibria  in future work. We also hope that this work may spark the interest of using the moving coordinates  among researchers.

We conclude this paper with a simple application of the results of stability to the  Sun-Jupiter system ($\mu_{SJ}\approx 9.538753\cdot 10^{-4}$) and  Earth-Moon system ($\mu_{EM}\approx 0.0121506$).
\begin{figure}
  \center
  \includegraphics[width=10cm]{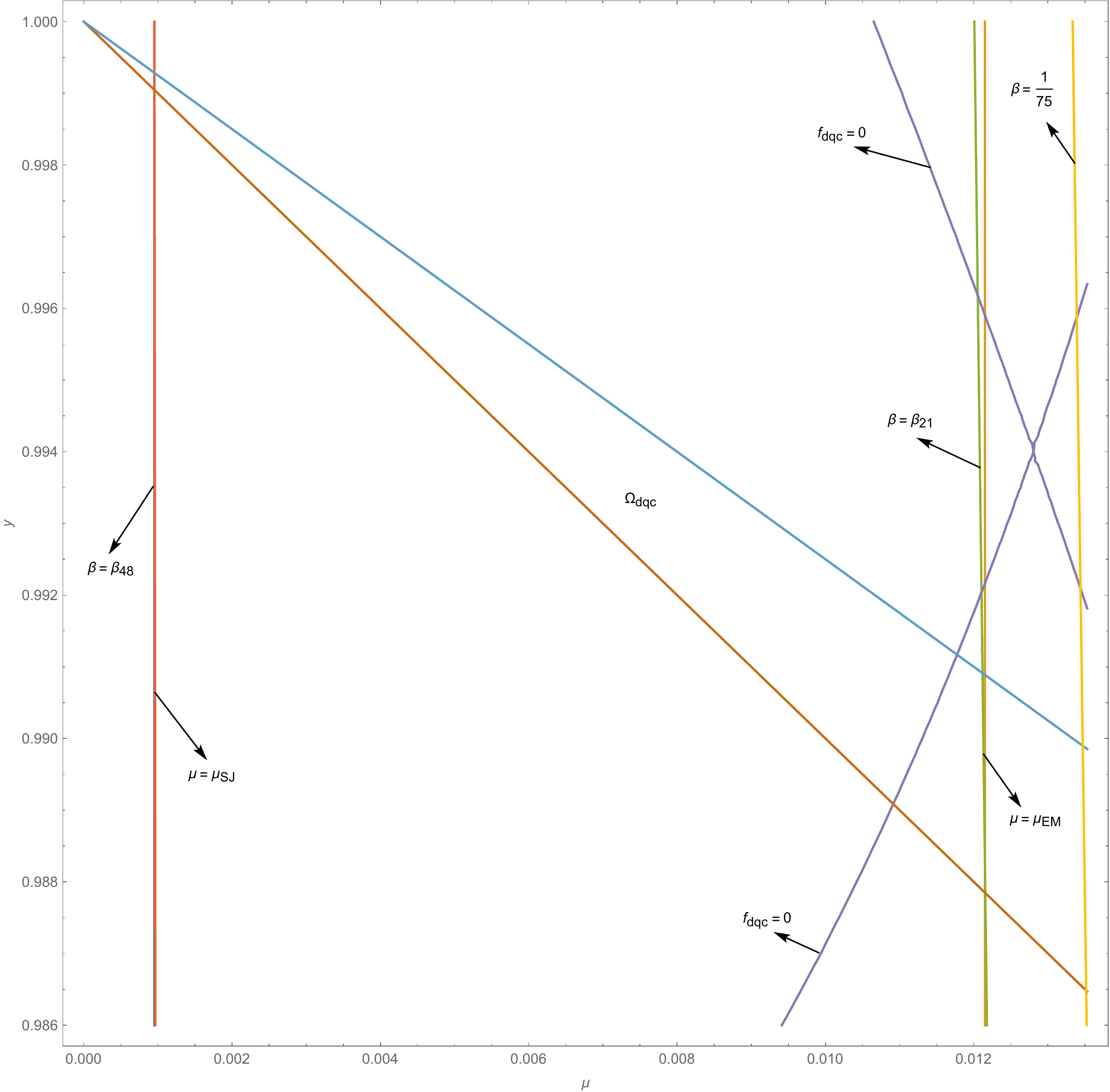}\\
  \caption{the  Sun-Jupiter system and Earth-Moon system}\label{Sun-JupiterEarth-Moon}
\end{figure}

First, each of 
Lagrange relative equilibria of the two systems is KAM stable, however we claim that the KAM stability on the  Sun-Jupiter system is stronger than on the Earth-Moon system in some sense. In fact, we can claim that the speed of  relative measure of  KAM invariant tori  tending to 1  for the  Sun-Jupiter system  is much faster than for the Earth-Moon system. Since the mass parameter $\beta$  first entering into  the stability
region,  for the Sun-Jupiter system, is $\beta_{48} \approx 9.530527\cdot 10^{-4}$ corresponding to a resonance
relation of order $48$;  and, for the Earth-Moon system, is $\beta_{21} \approx 0.0120078$ corresponding to a resonance
relation of order $21$.

For effective  stability, we believe that  each of
Lagrange relative equilibria of the two systems is exponentially stable, although the  Birkhoff normal
form $\mathcal{H}_{B4}$ of degree $4$ is  directionally quasi-convex for the  Sun-Jupiter system but not for the Earth-Moon system. Indeed, we believe that one can prove the exponential  stability of the Earth-Moon system by further calculating  its Birkhoff normal
form $\mathcal{H}_{B6}$ of degree $6$. On the other hand, we also believe that the stability  worsens if the condition of directional quasi-convexity for $\mathcal{H}_{B4}$ is violated, that is, the exponential stability on the  Sun-Jupiter system should be stronger than on the Earth-Moon system in some sense. As a matter of fact, in his
celebrated 1977 article \cite{nekhorshev1977exponential}, Nekhoroshev also conjectured that different steepness properties should lead to numerically observable differences in the stability times, although it is not easy to prove this.

All in all, the work makes us believe that the stability of Lagrange relative equilibrium for the Earth-Moon system is weaker than for the Sun-Jupiter system. This should be partially one of the reasons that  nobodies have been ever observed to gravitate
around Lagrange triangular  of the Earth-Moon system, but there are the well known Trojan asteroids  around Lagrange triangular  of the Sun-Jupiter system.

\indent\par


\section*{Appendix: on Lagrangian Dynamical Systems}\label{LagrangianDynamicalsystems}
      \indent\par
For the sake of readability, we sketchily give the theory of Lagrangian dynamical systems used in deducing the general equations of motion. The exposition follows \cite{zbMATH03601149,zbMATH05031968}, to which we refer the reader for proofs and details.

\begin{definition}
Let $\mathcal{M}$ be a differentiable manifold, $\mathrm{T}\mathcal{M}$ its tangent bundle, and $\mathcal{L}: \mathrm{T}\mathcal{M} \rightarrow \mathbb{R}$
a differentiable function. A map $\gamma: \mathbb{R} \rightarrow \mathcal{M}$ is called a motion in the Lagrangian
system with configuration manifold $\mathcal{M}$ and Lagrangian function $\mathcal{L}$ if $\gamma$ is an
extremal of the Lagrangian action functional
\begin{equation}
\mathcal{A}(\gamma) = \int^{t_2}_{t_1}{ \mathcal{L}(\dot{\gamma}(t)) dt}.\nonumber
\end{equation}
where $\dot{\gamma}$ is the velocity vector $\dot{\gamma}(t) \in \mathrm{T}_{\gamma(t)}\mathcal{M}$.
\end{definition}

\begin{theorem}\label{Euler-Lagrange
equations}
The evolution of the local coordinates $q = (q_1, \cdots, q_n)$ of a point $\gamma(t)$
under motion in a Lagrangian system on a manifold satisfies the \emph{Euler-Lagrange
equations}
\begin{equation}
\frac{d}{dt}\frac{\partial \mathcal{L}}{\partial \dot{q}} = \frac{\partial \mathcal{L}}{\partial {q}} ,\nonumber
\end{equation}
where $\mathcal{L}(q, \dot{q})$ is the expression for the function $\mathcal{L}: \mathrm{T}\mathcal{M} \rightarrow \mathbb{R}$ in the coordinates
$q$ and $\dot{q}$ on $\mathrm{T}\mathcal{M}$.
\end{theorem}

Theorem \ref{Euler-Lagrange
equations} yields a quick method
for writing equations of motion in various coordinate systems, even in larger class of coordinate transformations which  contain time. Indeed, to write the equations of motion in a new coordinate system, it is sufficient to
express the Lagrangian function in the new coordinates. In fact, we have
\begin{theorem}
If the orbit $\gamma: q = \varphi(t)$ of Euler-Lagrange equations $\frac{d}{dt}\frac{\partial \mathcal{L}}{\partial \dot{q}} = \frac{\partial \mathcal{L}}{\partial {q}}$ is written as $\gamma: Q = \Phi(t)$ in the local coordinates $Q, t$ (where $Q =
Q(q, t)$), then the function $\Phi(t)$ satisfies Euler-Lagrange equations $\frac{d}{dt}\frac{\partial \tilde{\mathcal{L}}}{\partial \dot{Q}} = \frac{\partial \tilde{\mathcal{L}}}{\partial {Q}}$, where $\tilde{\mathcal{L}}(Q, \dot{Q}, t) = \mathcal{L}(q,\dot{q}, t)$.

\end{theorem}

\begin{remark}
By the additional dependence of the
Lagrangian function on time:
\begin{equation}
\mathcal{L}: \mathrm{T}\mathcal{M}\times \mathbb{R} \rightarrow \mathbb{R}   ~~~~~~\mathcal{L}=\mathcal{L}(q,\dot{q},t),\nonumber
\end{equation}
one can consider a Lagrangian nonautonomous system and the results above are also valid.
  \end{remark}

\begin{definition}
In mechanics,  $\frac{\partial \mathcal{L}}{\partial \dot{q}}$ are called generalized
momenta, $\frac{\partial \mathcal{L}}{\partial {q}}$ are called generalized forces.
\end{definition}

\begin{definition}
Given a Lagrangian function $\mathcal{L}(q,\dot{q},t)$,  a coordinate $q_j$ is called ignorable (or cyclic ) if it does not enter into the
Lagrangian:  $\frac{\partial \mathcal{L}}{\partial q_j}=0$.
\end{definition}

\begin{theorem}
The generalized momentum corresponding to an ignorable coordinate is
conserved: $p_j = \frac{\partial \mathcal{L}}{\partial \dot{q}_j}=const$.
\end{theorem}

The following  content is \emph{Routh's method for eliminating
ignorable coordinates}.\\
Suppose that the Lagrangian $\mathcal{L}(q,\dot{q},\dot{\xi})$ does not involve the coordinate $\xi$, i.e., $\xi$ is ignorable. Using the equality $\frac{\partial \mathcal{L}}{\partial \dot{\xi}}=c$ we represent the  velocity
$\dot{\xi}$ as a function of $q,\dot{q}$ and $c$. Following Routh we introduce the function
\begin{equation}
\mathcal{R}_c(q,\dot{q})=\mathcal{L}(q,\dot{q},\dot{\xi}) - c \dot{\xi}|_{q,\dot{q},c} (= \mathcal{L}(q,\dot{q},\dot{\xi}) - \frac{\partial \mathcal{L}}{\partial \dot{\xi}}  \dot{\xi}|_{q,\dot{q},c} ).\nonumber
\end{equation}

\begin{theorem}
A vector-function $(q(t),\xi(t))$ with the constant value of generalized momentum $\frac{\partial \mathcal{L}}{\partial \dot{\xi}}=c$ satisfies the Euler-Lagrange
equations $
\frac{d}{dt}\frac{\partial \mathcal{L}}{\partial \dot{q}} = \frac{\partial \mathcal{L}}{\partial {q}}$
 if and only
if $q(t)$ satisfies
the Euler-Lagrange
equations $\frac{d}{dt}\frac{\partial \mathcal{R}_c}{\partial \dot{q}} = \frac{\partial \mathcal{R}_c}{\partial {q}}. $
\end{theorem}

{\bf Example. } Compare Newton's equations (\ref{eq:Newton's equation1})
\begin{equation*}
\frac{d{m_k \dot{\mathbf{r}}_k}}{dt} - \frac{\partial \mathcal{U}}{\partial \mathbf{r}_k}=0,
\end{equation*}
with the Euler-Lagrange equations
\begin{equation}
\frac{d}{dt}\frac{\partial \mathcal{L}}{\partial \dot{q}} - \frac{\partial \mathcal{L}}{\partial {q}}=0 \nonumber
\end{equation}
of  Lagrangian
system with configuration manifold $\mathcal{X}$ and Lagrangian function
\begin{displaymath}
\mathcal{L}  =  \mathcal{K} + \mathcal{U}.
\end{displaymath}
\begin{theorem}
Motions of the mechanical system (\ref{eq:Newton's equation1}) coincide with extremals of
the functional $\mathcal{A}(\mathbf{r})= \int^{t_2}_{t_1}{ \mathcal{L}dt}.$

\end{theorem}


\newpage

\bibliographystyle{plain}



\end{document}